\documentclass{article}
\usepackage{shapepar,bm}
\usepackage[dvips,arrow,matrix,ps,line]{xy}
\usepackage{theorem,amsfonts,amssymb,amscd, amsmath}
\usepackage{geometry}
\usepackage{ytableau,youngtab}

\usepackage{makeidx}
\usepackage{stmaryrd,setspace,palatino,upgreek}

\usepackage{graphicx}
\usepackage{graphics}

\usepackage{epstopdf}

\usepackage{tikz}
\usetikzlibrary{
  decorations.pathreplacing,
  arrows.meta,
  positioning
}
\usepackage{amssymb,graphics,hyperref}

\hypersetup{
    colorlinks = true,
    linkcolor = blue,
    anchorcolor = red,
   citecolor = blue,
    filecolor = red,
    pagecolor = red,
    urlcolor = blue}

\newcommand{\mathsym}[1]{{}}

\makeindex

\def\bfd{{\mathbf d}}

\def\bbS{\mathbb S}
\newcommand{\llb}{\llbracket}
\newcommand{\rrb}{\rrbracket}
\def\bmtx{\begin{matrix}}
\def\emtx{\end{matrix}}
\def\NN{\mathbb N}
\def\id{\mathrm{id}}

\def\bfe{{\mathbf e}}
\def\bfz{{\mathbf z}}

\def\ovsig{\overline{\sigma}}

\def\bwV{{\bigwedge\hskip-3.5pt V}}

\def\bfX{{\mathbf X}}

\def\bfw{{\mathbf w}}
\def\Qcal{\mathcal Q}

\def\d{\partial}
\def\bmd{\bm\d}

\def\Gcal{\mathcal G}

\def\bfx{{\mathbf x}}
\def\bfX{{\mathbf X}}

\def\bff{{\mathbf f}}

\def\ZZ{\mathbb Z}

\def\CC{\mathbb C}

\def\bff{{\bf f}}

\def\Scal{{\mathcal S}}

\def\QQ{\mathbb Q}

\def\PP{\mathbb P}

\def\cocoa{{\hbox{\rm C\kern-.13em o\kern-.07em C\kern-.13em o\kern-.15em A}}}

\def\Dcal{\mathcal D}

\def\Span{{\rm Span}}

\def\bft{{\bf t}}

\def\ovD{\overline{D}}
\def\ovDcal{\overline{\Dcal}}

\def\End{\mathrm{End}}

\def\blamb{{\bm \lambda}}
\def\bmu{{\bm \mu}}
\def\bnu{{\bm \nu}}
\def\Pcal{{\mathcal P}}

\def\ep{{\epsilon}}

\def\w2M{\bigwedge^2M}

\def\w{\wedge }
\def\bw{\bigwedge }

\protect
\protect
\protect
\protect

\protect
\protect

\def\sra{\rightarrow}
\def\lra{\longrightarrow}

\def\proof{\noindent{\bf Proof.}\,\,}
\def\qed{{\hfill\vrule height4pt width4pt depth0pt}\medskip}
\def\be{\begin{equation}}
\def\ee{\end{equation}}
\def\bclm{\begin{claim}}
\def\eclm{\end{claim}}
\def\beqn{\begin{eqnarray}}
\def\eeqn{\end{eqnarray}}
\def\beqn*{\begin{eqnarray*}}
\def\eeqn*{\end{eqnarray*}}

\theoremstyle{change}
\theorembodyfont{\rmfamily}

\newtheorem{claim}{}[section]

\global
\global

\def\no@breaks#1{{\def\\{ \ignorespaces}#1}}    

  \makeatletter
\def\cleardoublepage{\clearpage\if@twoside \ifodd\c@page\else
\hbox{} \thispagestyle{empty}
\newpage
\if@twocolumn\hbox{}\newpage\fi\fi\fi} \makeatother

\usepackage{eso-pic,graphicx}
\makeatletter
\newcommand\BackgroundPicture[2]{%
  \setlength{\unitlength}{1pt}%
  default \put(0,\strip@pt\paperheight){%
  \parbox[t][\paperheight]{\paperwidth}{%
    \vfill
     \centering \includegraphics[angle=#2, width=15cm, height=15cm,  bb=0 0 150 150]{#1}
    \vfill
}}} %
\makeatother

\providecommand{\bysame}{\leavevmode\hbox to3em{\hrulefill}\thinspace}

\def\bfX{{\mathbf X}}
\title{On Two Algebraic Realizations of Schubert Calculus}
\author{André Contiero, Letterio Gatto, Parham Salehyan}

\begin{document}

\maketitle
\begin{abstract}
\noindent Schubert calculus on complex Grassmannians can be played by means of differential operators acting on Schur polynomials or Vertex Operators acting on exterior algebras. In this paper we develop this point of view systematically and complement it with a parallel exterior-algebra formalism, leading to what we call the \emph{bosonic} and the \emph{fermionic Schubert calculi}.
\noindent
The two alluded realizations are related by (a finite type version of) the boson--fermion correspondence, thereby providing a unified framework connecting Schubert calculus and integrals on the Grassmannian, symmetric functions, exterior algebras and the representation theory of symmetric groups.
\end{abstract}

\setstretch{1.25}

\section{Introduction}

On the one hand, consider the $\QQ$-algebra
$
B:=\QQ[\bfx]$ in the indeterminates  $\bfx=(x_1,x_2,\ldots)$,
the polynomial ring in infinitely many indeterminates, endowed with  its  Schur $\QQ$-basis $(\Scal_\blamb(\bfx)_{\blamb\in\Pcal}$ parametrized by the set $\Pcal$ of all the integer partitions (see Section~\ref{sec1:parti}). It carries a natural action of the (non-commutative) infinite Weyl algebra
\[
A_\infty(\QQ)
=
\displaystyle{\left\langle
x_i,\frac{\partial}{\partial x_j}
\right\rangle_{i,j\geq 1}},
\] 
generated by the \emph{creation operators} (multiplication by $x_i$) and the \emph{annihilation operators} (partial derivatives $\partial/\partial x_i$), subject to the commutation relations
$
\left[{\partial}/{\partial x_i},x_j\right]=\delta_{ij}$,
$[x_i,x_j]=0$,
\linebreak $\left[{\partial}/{\partial x_i},{\partial}/{\partial x_j}\right]=0.
$

On the other hand, let $G(r,n)$ denote the complex Grassmannian parametrizing $r$-dimensional 
subspaces of $\CC^n$. It is a smooth complex projective variety whose intersection theory is  governed by the classical {\em Schubert calculus}. The latter describes the action of the cohomology ring
$
H^*(G(r,n),\QQ)
$
on the homology
$
H_*(G(r,n),\QQ)
$ via the {\em cap product}. The cohomology ring is generated by the special Schubert classes
$
\sigma_i=c_i(\Qcal_r),
$
the Chern classes of the universal quotient bundle, whereas the homology possesses a distinguished basis
$
\{\Omega^\blamb\},
$
indexed by the set $\Pcal_{r,n}$ of the partitions whose Young diagrams are contained in an $r\times(n-r)$ rectangle.  

The singular cohomology of the Grassmannian glues the first actor (the $\QQ$-algebra $B$) with 
a third  main character of this article: the exterior algebra $\bw V$ of the vector space 
$V=\QQ[X]$ of the polynomials in one indeterminate,  acted on by a certain Lie sub-algebra 
$gl(V)$ of the Lie algebra of all the endomorphisms of $V$ (Sect.~\ref{sec:sec49}). The 
purpose of this paper is to discuss the bridges among these three realms, in connection to 
some classical subjects like the representation theory of the Heisenberg algebra and that of 
the symmetric group on $d$ letters.

The Ariadne thread we follow, indeed our original motivation, is the notion of integration on a Grassmannian. By this, we mean the numerical invariants obtained by evaluating top-dimensional intersections such as
\[
\int_{G(r,n)}
\sigma_1^{i_1}\cdots\sigma_k^{i_k}\cap [G(r,n)].
\]
These numbers vanish unless
$
i_1+2i_2+\cdots+ki_k=r(n-r),
$
and play a central role in enumerative geometry. Historically, intersection numbers of this type provided one of the paradigms which later inspired the study of intersection theory on moduli spaces of stable curves and tautological classes as explicitly declared, e.g. in \cite{MumfEnum} and effectively implemented in \cite{Wit2dim}.

There is a vast literature concernig Schubert Calculus: the foundational part is easily available in many excellent textbooks, such as the classical \cite{Ful,GH}, the book \cite{Manivelsym} for a more symmetric functions oriented reading and the milestone paper by Kleiman and Laksov \cite{KleLak1}.

Inspired by previous contributions to the subject within the same spirit, such as in \cite{SCHSD,SCGA,gln,BeCoGaVi,GSCH}, and \cite{GaRow1} in the more general framework of semi-algebras, the purpose of this paper is to revisit Schubert calculus from two complementary viewpoints, which we shall call the \emph{bosonic Schubert calculus} and the \emph{fermionic Schubert calculus}. The terminology is motivated by the fact that the ring $B$  carries the standard bosonic action of a Weyl algebra through creation and annihilation operators, whereas the exterior algebra carries the corresponding fermionic action through {\em Schubert derivations} (as in~\cite{SDIWP}). Although the two formalisms are ultimately equivalent, they emphasize different aspects of the theory.

The bosonic Schubert calculus is developed on the polynomial ring
$
B=\QQ[x_1,x_2,\ldots],
$ (a representation space for many relevant infinite dimensional Lie algebras),
endowed with its Schur basis
$
\{\Scal_\lambda(\bfx)\}_{\lambda\in\Pcal}.
$, parametrized by the set $\Pcal$ of all of the integer partitions .
We regard $B$ as a universal receptacle containing the homologies of all Grassmannians at once. More precisely, for each pair $(r,n)$ there is a natural projection
$
\pi_{r,n}:B\longrightarrow H_*(G(r,n),\QQ),
$
sending the Schur polynomial $\Scal_\blamb(\bfx)$ to the corresponding Schubert cycle.

A key role throughout the paper is played by a simple (and known!) determinantal identity, stated in Lemma~\ref{lem:giapie}. It provides a common source for the Giambelli- and Jacobi--Trudi-type formulas appearing in the paper and yields Pieri-type identities in a uniform way. In this sense, several classical formulas of Schubert calculus and symmetric function theory arise from one elementary mechanism, in our opinion not very exploited in literature.

The ``bosonic'' guiding principle is that intersection theory on Grassmannians admits a realization through the natural action on $B$ of the infinite Weyl algebra recalled at the beginning of this introduction, where the multiplication operators $p\mapsto x_ip$ play the role of creation operators, whereas the partial derivatives
$
p\mapsto \displaystyle{\frac{\partial p}{\partial x_i}}
$
play the role of annihilation operators. 

Within this framework, Schubert cohomology classes become explicit polynomial expressions in the operators $\partial/\partial x_i$ (Cf. Section~\ref{sec:sec313}). Thus the partial derivatives may be regarded as elementary cohomological operators, while the usual Schubert classes appear as polynomial expressions in them.
More precisely, the total Chern class of the universal quotient bundle is represented by the differential operator
$
\displaystyle{\exp\left(
\sum_{i\geq 1}
\frac{t^i}{i}
\frac{\partial}{\partial x_i}
\right)}
$,
which is  the ``differential half'' of the standard vertex operator occurring in the boson--fermion correspondence (see e.g. \cite[Proposition 5.1]{KaRai}). Equivalently, if
\[
\sum_{i\geq 0}
\Scal_i(\widetilde{\partial})t^i
=
\exp\left(
\sum_{i\geq 1}
\frac{t^i}{i}
\frac{\partial}{\partial x_i}
\right),
\]
then the action of the special Schubert class $\sigma_i$ on the homology vector space of the Grassmannian is modeled by the action of the operators $\Scal_i(\widetilde{\partial})$ on $B$.

Thus intersection integrals become evaluations of differential operators. For instance, if $\sum_{j=1}^kji_j=r(n-r)$
\[
\int_{G(r,n)}
\sigma_1^{i_1}\cdots\sigma_k^{i_k}
\cap [G(r,n)]
=
\Scal_1(\widetilde{\partial})^{i_1}
\cdots
\Scal_k(\widetilde{\partial})^{i_k}
\Scal_{((n-r)^r)}(\bfx).
\]

In particular (Corollary \ref{cor:cor31}), the Pl\"ucker degree of the Schubert variety indexed by a partition $\blamb$ is given by 
\[
\omega_\blamb
=
\frac{\partial^{|\blamb|}}
{\partial x_1^{|\blamb|}}
\Scal_\blamb(\bfx),
\]
a formula which  makes explicit the relation between degrees of Schubert varieties, Schur polynomials and the number of standard Young tableaux of a fixed shape \cite[Ch.~4]{Fulyoung}.

The {\em fermionic Schubert calculus}, instead,  is developed on the exterior algebra of the vector space
$
V=\QQ[X].
$
A central role in the fermionic picture is played by the {\em vacuum vector}
\[
\bfX^r(0)
=
X^{r-1}\wedge X^{r-2}\wedge\cdots\wedge X^0,
\]
whose geometric interpretation is that of the fundamental class of the infinite dimensional Grassmannian $G(r,\infty)$ and provides the bridge between the bosonic and fermionic formalisms. In particular, the evaluation of polynomials of $B$ against the vacuum, plays the role of intersection with the fundamental class.

The exterior power $\bigwedge^rV$ possesses a double interpretation. On the one hand, through the distinguished isomorphism
$
B_r
\longrightarrow
\bigwedge^rV,
$
 mapping $
\Scal_\lambda(\bfx_r)
\longmapsto
\bfX^r(\blamb),
$
it may be regarded as a copy of the cohomology ring
$
B_r\cong H^*(G(r,\infty),\QQ).
$
by transporting the ring structure.

On the other hand, the same space is naturally a free $B_r$-module generated by the vacuum vector
$
\bfX^r(0),
$
and this module structure models the cap product of cohomology with homology. In this sense $\bigwedge^rV$ may be viewed simultaneously as a fermionic model for cohomology and as a homology module acted on by cohomology.

Within this framework Schur polynomials are replaced by wedge monomials
\[
\bfX^r(\blamb)
=
X^{r-1+\lambda_1}
\wedge
X^{r-2+\lambda_2}
\wedge
\cdots
\wedge
X^{\lambda_r},
\]
and cohomological operations are realized through Schubert derivations, namely distinguished {\em Hasse--Schmidt derivations} \cite{SCHSD,HSDGA} acting on exterior powers. The boson--fermion correspondence identifies the bosonic and fermionic Schubert calculi.

The exterior-algebra formalism also provides finite-rank replacements for the partial derivatives. Namely, the unique derivations
$
\delta(X^{-i})
$
of the exterior algebra such that $X^{-i}X^j=X^{j-i}$ if $j\geq i$ and $0$ otherwise, 
act on exterior powers as annihilation operators and turn, in the stable limit, into the partial derivatives
$
\frac{\partial}{\partial x_i}
$ acting on $B$.
At finite rank, their commutation relations with the creation operators $\delta(X^j)$ contain explicit {\em boundary terms corrections} explaining why the usual Heisenberg relations fail on finite exterior powers and why and how they reappear  in the stable limit.

One nice output we have not been able to find in any previous literature, probably due to the fact of being still unusual playing Schubert Calculus on exterior algebras,  is that the same annihilation operators recover the irreducible characters of symmetric groups. More precisely (Theorem~\ref{thm:thm431}), if
$
\bmu=(1^{m_1}2^{m_2}\cdots k^{m_k})\vdash d
$
and $\blamb\vdash d$, then
\[
\delta(X^{-1})^{m_1}
\delta(X^{-2})^{m_2}
\cdots \delta(X^{-k})^{m_k}
\bfX^r(\blamb)
=
\chi^\blamb_\bmu\,\bfX^r(0),
\]
for every $r\geq \ell(\blamb)$. Thus character values arise as {\em vacuum coefficients} of products of fermionic annihilation operators (basically partial derivatives). Equivalently, they may be interpreted as fermionic integrals, in the same sense that Pl\"ucker degrees arise as ordinary intersection numbers. From this perspective, Schubert degrees and irreducible character values appear as bosonic and fermionic manifestations of the same integration principle. In passing we obsere that as a further direction of research, the set of partitions offer an entry point to define some {\em Littlewood-Richardson} semi-rings (basically the monoid of partitions equipped with the prdocut of partitions having the Littlewood-Richardson coefficients as structural constants), which can be likely analyzed within the spectacles offered by the theory of systems by L.~Rowen (see e.g. \cite{RoInfo}, \cite{RoSemiPa}, \cite{RoNega}), by enhancing the and making more effective the philosophy already proposed in \cite{GaRow1}.

The paper is organized as follows. In Section~\ref{sec1:preli} we recall the necessary 
background on partitions, Schur polynomials and exterior powers, emphasizing the determinantal 
identity which will be used throughout. In Section~\ref{sec3:diff} we develop the bosonic 
Schubert calculus and describe the action of Schubert classes by differential operators. We 
also discuss some consequences for Pl\"ucker degrees of Schubert varieties and related 
identities in symmetric function theory. In Section~\ref{sec:sec4}  the {\em fermionic} 
version of Schubert calculus is introduces through Schubert derivations on exterior algebras. 
Certain  boundary Heisenberg relationsare discussed as well and finally relate annihilation operators to irreducible characters of symmetric groups (Theorem~\ref{thm:thm431}).

\medskip
\noindent
{\bf Acknowledgment.} The second named author (L.G.) is grateful to Inna Scherbak for her unreserved support along many years, inclusive in encouraging part of the present project,  Renato Vidal Martins for many inspiring discussions and Edoardo Ballico for generous advises.

\section{Preliminaries}\label{sec1:preli}
\claim{\bf Partition.} \label{sec1:parti} A partition is a non increasing monotonic sequence of non-negative integers \linebreak
$
\blamb=(\lambda_1\geq \lambda_2\geq \cdots \geq \lambda_r\geq 0)
$
such that $\lambda_i=0$ for all but finitely many $i\in\NN$.  Let $\Pcal$ be the additive 
monoid of all the partitions. If $\blamb\in\Pcal$, its length $\ell(\blamb)$ is the 
cardinality of its support, whose elements are called {\em parts}: $
\ell(\blamb):=\sharp\{i\mid \lambda_i\neq 0\},
$, while $|\blamb|:=\sum_{i\geq 1}\lambda_i$ is its weight. For all integer $i\geq 0$, one 
denotes by $m_i(\blamb):=\{j\,|\, \lambda_j=i\}$ so that a partition $\blamb$ can be 
equivalently written as a sequence $(1^{m_1}2^{m_2}\cdots)$ where all $m_i:=m_i(\blamb)$ are 
zero but finitely many, which are omitted from the sequence. For instance 
$(5,2,2,1)=(1^12^25^1)$. 
\claim{} \label{sec:sec32} One denotes by $\Pcal_r$ the sub-monoid of $\Pcal$ whose elements are the partitions 
of length at most $r$: $\Pcal_r:=(\blamb\in\Pcal\,|\,\ell(\blamb)\leq r\}\subseteq \NN^r$. By 
$\Pcal_{r,n}$ one understands the set of all $\blamb\in\Pcal_r$ such that $n-r\geq \lambda_1$. They can be 
characterized by the property that their Young diagram is contained in a $r\times (n-r)$ 
rectangle.
\claim{} If $\blamb\in\Pcal_r$, its Young diagram is an array of $r$ left-aligned rows of boxes, where the $i$-th row contains $\lambda_i$ boxes.
A {\em Young tableau} of shape $\blamb$ is a filling of the boxes of the Young diagram with 
integers. It is called a \emph{standard Young tableau} if the entries are the integers $1,2,
\dots,|\lambda|$, each appearing exactly once, and increasing along each row (from left to 
right) and along each column (from top to bottom). The important combinatorial and 
representation theoretical interpretation of the number of Young tableaux of shape $\blamb$ 
is that it coincides with the dimension of the {\em Specht module} $\bbS^\blamb$, which is the 
irreducible representation of the symmetric group on $d=|\blamb|$ elements associated to the 
partition $\blamb$ (see \cite{Fulyoung}). Such a dimension can be computed by using the 
celebrated {\em hook lenght fromula} by Frame-Robinson and Thrall  \cite{thrall}. We will 
look at the dimension of $\bbS^\blamb$ rather as the degree $\omega_\blamb$ of the Schubert variety in a Grassmannian $G(r,n)$ ($r\geq \ell(\blamb)$) of 
dimension $|\blamb|$ associated to the partition $\blamb$.

If $\blamb\in \Pcal_{r,n}$, we denote by $\blamb^c$ the partition whose Young diagram is the 
complement of the Young diagram of $\blamb$ in a $r\times (n-r)$ rectangle.
For example, the complement of $\blamb=(3,3,2,1)$ in the $4\times 3$ rectangle is $
\blamb_1^c(2,1)$, whereas its complement in the $5\times 4$ rectangle is $\blamb_2^c=(4,3,2,
1,1)$. The conjugated $\blamb'$ of $\blamb\in\Pcal_r$ is the partition whose Young diagram is 
the tranposed of $Y(\blamb)$.
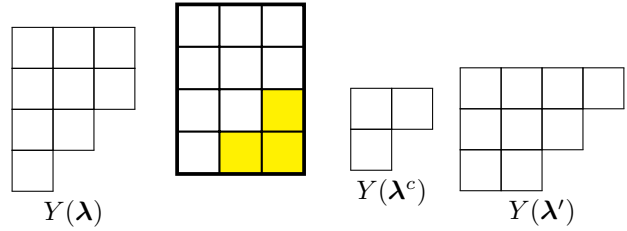
\[
\begin{matrix}
\bmtx\begin{ydiagram}{3,3,2,1}\end{ydiagram} \cr Y(\blamb)\cr \cr\cr\emtx\,\,\,& 
\begin{tikzpicture}[scale=0.56]
\def\s{1}

\fill[yellow] (2*\s,-3*\s) rectangle ++(\s,\s);
\fill[yellow] (1*\s,-4*\s) rectangle ++(\s,\s);
\fill[yellow] (2*\s,-4*\s) rectangle ++(\s,\s);

\draw[step=\s, line width=0.8pt] (0,0) grid (3*\s,-4*\s);

\draw[line width=1.4pt] (0,0) rectangle (3*\s,-4*\s);
\end{tikzpicture}&\,\,\,\,\bmtx\ydiagram{2,1}\cr Y(\blamb^c)\cr\cr\emtx&\bmtx\begin{ydiagram}{4,3,2}\end{ydiagram}\cr Y(\blamb')\cr\cr\emtx
%
%
%
\end{matrix}
\]

\claim{\bf Symmetric polynomials.} Let $\Lambda_r:=\QQ[z_1,\ldots,z_r]^{S_r}$ be the ring of symmetric polynomials in the $r$ indeterminates $\bfz_r:=(z_1,\ldots,z_r)$.  By the fundamental theorem of symmetric functions, $\Lambda_r$ is the polynomial $\QQ$-algebra generated by the elementary symmetric functions
$
e_1(\bfz_r),\ldots,e_r(\bfz_r)
$, 
defined by $\sum_{j\geq 0}(-1)^je_j(\bfz_r)t^j=\prod_{i=1}^r(t-z_i)$. A typical 
basis element of $\Lambda_r$ is then $e_1^{i_1}(\bfz_r),\cdots, e_r^{i_r}(\bfz_r)$. The 
partition $(1^{i_1},\ldots,r^{i_r})$ is conjugated to a partition of length at most $r$, then 
$\Lambda_r$ has a basis parametrized by partitions. 

For all $s\geq r$, there is an obvious map  $\Lambda_s\sra \Lambda_r$, defined by putting $z_j=0$ if $j>r$. The ring of symmetric functions (in infinitely many indeterminates) is the projective limit $\Lambda=\underset{\to}{\Lambda_r}$ in the category of graded algebras (as in \cite[Remark 1, p.~19]{MacDonald}). 

One distinguished basis is provided by 
the {\em Schur polynomials} $s_\blamb(\bfz_r)$, which one defines as the ratio of two skew symmetric polynomials, 
namely
\be
s_\blamb(\bfx)={\det(z_j^{i-1+\lambda_{r-i+1}})\over \Delta_0(\bfz_r)}\label{eq2:Vand}
\ee
where $\Delta_0(\bfz_r):=\prod_{1\leq i<j\leq r}(z_i-z_j)$ is the {\em Vandermonde} determinant.

A standard argument shows that the polynomials $s_\blamb(\bfz_r)$ are linearly independent, 
and  hence they do form a basis of $\Lambda_r$. The following statement is well known and is the combinatorial essence of all the Giambelli's and Pieri's type formulas in the literature. It will be repeatedly used along the exposition.

\bclm{\bf Lemma.}\label{lem:giapie} {\em 
Let \(f(z)=\sum_{m\in\mathbb Z}f_mz^m\) be a formal Laurent series with coefficients in a
commutative \(\QQ\)-algebra \(A\). For \(r\geq 1\), define \(\Delta_\blamb(\bff)\in A\) via the the equality
\[
\sum_{\blamb\in\Pcal_r}
\Delta_\blamb(\bff)\,s_\blamb(z_1,\ldots,z_r)=\bff(\bfz_r):=f(z_1)\cdots f(z_r)
\]
\noindent
Then
\be
\Delta_\blamb(f)=\det(f_{\lambda_j-j+i})=\begin{vmatrix}
f_{\lambda_1}&f_{\lambda_2-1}&\cdots&f_{\lambda_r-r+1}\cr
f_{\lambda_1+1}&f_{\lambda_2}&\cdots&f_{\lambda_r+r-2}\cr
\vdots&\vdots&\ddots&\vdots\cr
f_{\lambda_1+r-1}&f_{\lambda_2+r-2}&\cdots&f_{\lambda_r}.
\end{vmatrix}\label{eq6:pregiam}
\ee}

\eclm
\proof
First of all one notices that multiplying both sides of the equality by the Vandermonde $\Delta_0(\bfz_r)=\prod_{i<j}(z_i-z_j)$ one obtains
$$
\begin{vmatrix}
f(z_1)&f(z_2)&\cdots&f(z_r)\cr
z_1f(z_1)&z_2f(z_2)&\cdots&z_rf(z_r)\cr
\vdots&\vdots&\ddots&\vdots\cr
z_1^{r-1}f(z_1)&z_2^{r-1}f(z_2)&\cdots&z_r^{r-1}f(z_r)
\end{vmatrix}=\sum \Delta_\blamb(\bff)\det(z_j^{i-1+\lambda_{r-i+1}})
$$
It follows that the coefficient of $z_1^{r-1+\lambda_1}\cdots z_r^{\lambda_r}$ is precisely $\Delta_\blamb(f)$, while on the left hand side is $\det(f_{\lambda_j-j+i})$, as desired.

\qed

\claim{} For all $i\geq 0$ and $\blamb\in\Pcal_r$, let
$$
PF_i(\blamb)=\{\bmu\in\Pcal_r\,|\, \mu_1\geq \lambda_1\geq 	\ldots\geq \mu_r\geq\lambda_r \,\,\mathrm{and}\,\, |\bmu|=|\blamb|+i\}
$$
and
$$
PF_{-i}(\blamb)=\{\bmu\in\Pcal_r\,|\, \lambda_1\geq \mu_1\geq 	\ldots\geq \lambda_r\geq\mu_r \,\,\mathrm{and}\,\, |\bmu|=|\blamb|-i\}
$$
\bclm{\bf Remark.}\label{rmk:remark39} It is important to notice that $\bmu\in PF_i(\blamb)$ if and only if $\blamb\in PF_{-i}(\bmu)$. Moreover it is worth to notice that if  $\blamb\in\Pcal_{r,n}$ (Cf.~\ref{sec:sec32}), then $PF_{-i}(\blamb)=PF_i(\blamb^c)$.\eclm

Let now $f_i$ and $\Delta_\blamb(\bff)$ as in Lemma~\ref{lem:giapie}. Then:
\bclm{\bf Proposition.}\label{propo:genpier} {\em   Pieri's formula holds:
\be 
f_i\Delta_\blamb(\bff)=\sum_{\bmu\in PF_i(\blamb)}\Delta_\bmu(\bff)\label{eq:fultonap}
\ee}
\eclm
\qed


\section{Schubert Calculus Through Differential Operators}
\label{sec3:diff}
The aim of this Section is twofold. First, one quickly recalls the 
classical Schubert Calculus backgorund for usual complex Grasmannians, 
essentially following \cite[p.~271]{Ful}. Next,  we show how the former is 
completely ruled by the structure of the ring $B$ as a module over its 
Weil Algebra $A_\infty$, generated by the multiplication by $x_i$ and the 
derivatives with respect $x_i$.
\claim{\bf Schubert Calculus, in short.} Let $G(r,n)$ be the complex 
Grassmannian parametrizing $r$-dimensional vector subspaces of $\CC^n$. It 
is a complex projective variety which has a celluar decomposition over the 
reals,  with affine cells in even dimension only \cite{GH,Ful}. According 
to  \cite[p.~271]{Ful}, one is lead to consider an arbitrary flag of 
vector subspaces of $\CC^n$
\[
F_\bullet(\blamb): \quad
0 \subseteq F_1 \subsetneq F_2 \subsetneq \cdots \subsetneq F_r \subseteq 
\CC^n
\]
such that
$
\dim F_i = i + \lambda_{r-i+1}.
$
Any such flag determine and is determined by a partition \linebreak $\blamb=(\lambda_1\geq\cdots\geq \lambda_r)\in 
\Pcal_{r,n}$ to which one attaches the Schubert variety associated to $F_\bullet$ 
of dimension $|\blamb|$:
\[
\Omega^\blamb(F_\bullet)
:=
\{\Lambda \in G(r,n)\mid \dim(\Lambda \cap F_i)\geq i\},
\]
which is the closure of a corresponding Schubert cell of the same 
dimension.
Its homology class in $H_*(G(r,n),\QQ)$ does not depend on the choice of $F_\bullet$ and it will be denoted by $\Omega^\blamb$. Each $\Omega_\blamb(F^\bullet)$ is a closed projective subvariety of $G(r,n)$ which in turn can be embedded à la Pl\"ucker in a projective space of dimension $N:={n\choose r}-1$. 
The other classical alternative descripton of the cellular decomposition of $G(r,n)$ is to consider a complete flag of vector subspaces of $\CC^n$
\[
\Gcal^\bullet: \CC^n = G^0 \supset G^1 \supset \cdots \supset G^n = 0
\]
such that $\operatorname{codim} G^i = i$.  
For $\blamb\in\Pcal_{r,n}$, let
\[
\Omega_\blamb(\Gcal^\bullet)
:=
\left\{
V \in G(r,n)\ \middle|\ 
\dim\bigl(V \cap G^{\,n-r+i-\lambda_i}\bigr) \ge i,\quad 1 \le i \le r
\right\}
\]
also called Schubert variety.
Its homology class $[\Omega_\blamb(\Gcal^\bullet)]\in H_*(G(r,n), \QQ)$ does not depend on the choice of the flag and is  very well known that
$$
H_*(G(r,n),\QQ):=\bigoplus_{\blamb\in \Pcal_r}\QQ\cdot \Omega^\blamb= \bigoplus_{\blamb\in \Pcal_r}\QQ\cdot\Omega_\blamb 
$$
as well as that 
\be 
\Omega^\blamb=\Omega_{\blamb^c},\label{eq:dualom}
\ee for $\blamb\in\Pcal_{r,n}$.

\claim{} The  Pl\"ucker degree $\omega_\blamb$ of the Schubert class $\Omega^\blamb$, geometrically speaking, is the number of reduced points of intersections with $|\blamb|$ hyperplanes in general position,  is independent of $n\geq r$. It coincides, e.g. by \cite[Theorem 2.31]{smirnov}, with the number of {\em standard Young tableaux} of shape $\blamb$ and we are going to learn how to compute it by a simple iteration of a partial derivative. 
\claim{} Singular homology is acted on by singular cohomology.
 Let
$
\Qcal_r \longrightarrow G(r,n)
$
be its universal quotient bundle. Denote by $\sigma_i=c_i(\Qcal_r)\in H^*(G(r,n),\QQ)$  the $i$-th Chern class of $\Qcal_r$, also said to be a {\em special Schubert cycle}. Let $c_t(\Qcal_r)=1+\sum_{i\geq 1}\sigma_it^i\in H^*(G(r,n),\QQ)[t]$ its Chern polynomial (of degree $n-r$). 
The classical Pieri formula  such as, e.g., \cite[p.~283]{GH} says that
\be
\sigma_i\cap \Omega_\blamb=\sum_{\bmu\in PF_i(\blamb)}\Omega_\bmu.
\ee
Due to \eqref{eq:dualom} one deduces
\bclm{\bf Proposition.}\label{prop:prop23} {\em The dual Pieri's rule holds:
\be
\sigma_i\cap \Omega^\blamb=\sum_{\bmu\in PF_{-i}(\blamb)}\Omega^\bmu.\label{eq5:dualp}
\ee
}
\eclm

\proof
It is based on the duality $\Omega_\blamb=\Omega^{\blamb^c}$, suitably applying  Remark~\ref{rmk:remark39}:
$$
\sigma_i\cap \Omega^\blamb=\sigma_i\cap \Omega_{\blamb^c}=\sum_{\bmu\in PF_i(\blamb^c)}\Omega_\bmu= \sum_{\bmu\in PF_i(\blamb^c)}\Omega^{\bmu^c}=\sum_{\bmu\in PF_{-i}(\blamb)}\Omega^\bmu.
$$

\qed


\medskip

\claim{}\label{sec:sec41} The purpose is now to illustrate how Schubert Calculus formalism can be phrased in terms of differential calculus. The main characters are the following related rings. They model, according to the context, the homology or the cohomology of the Grassmannians. 

For $r\in\NN\cup\{\infty\}$, let $B_r:=\QQ[e_1,\ldots,e_r]$, a polynomial ring in some ``elementary'' indeteminates $(e_1,\ldots,e_r)$. The convention is that $B_0=\QQ$. The ring $B_r$ comes endowed with some generic monic  polynomial $E_r(z)=1-e_1z+\cdots+(-1)^re_rz^r\in B_r[z]$ (for finite $r$), through which one defines a sequence $\bfx_r:=(x_1,x_2,\ldots)$ of $\QQ$-algebra generators, imposing the equality
$$
\sum_{i\geq 1}x_iz^i=1-E_r(z)+{(1-E_r(z))^2\over 2}+{(1-E_r(z))^3\over 3}+\cdots=-\log(1-E_r(z))\in B_r\llb z\rrb
$$
It follows that $\bfx_r:=(x_1,x_2,\ldots)$ generate  $B_r$ as a $\QQ$-algebra and we shall write, abusing notation, $B_r=\QQ[x_1,\ldots,x_r]$, because only  the first $r$ are algebraically independent: each $x_i$ is a homogeneous polynomial of weighted degree $i$ in $(x_1,x_2,\ldots,x_r)$.

\claim{} \label{sec:secB}

The equality
$$
\sum_{i\geq 0}\Scal_i(\bfx_r)z^i:=\exp(\sum_{i\geq 1}x_iz^i)
$$
defines {\em Schur} polynomials $(\Scal_i(\bfx_r))_{i\geq 0}$.
Since $\Scal_i(\bfx):=\displaystyle{x^i\over i!}+ g(\bfx_r)$, where $\deg_{x_1}g<i$, it follows that $\Scal_1(\bfx_r),\ldots, \Scal_r(\bfx_r)$ also generate $B_r$ as a $\QQ$-algebra. The quotient
$$
B_{r,n}:={B_r\over (\Scal_{n-r+1}(\bfx_r),\ldots, \Scal_n(\bfx_r))}\cong \bigoplus_{\blamb\in\Pcal_{r,n}}\QQ\Scal_\blamb(\bfx_r)
$$
is classically well known as the  presentation of the cohomology of the Grassmannian $G(r,n)$ where $\Scal_j(\bfx_r)$ play the role of the Chern classes of the universal quotient bundle of rank $n-r$.

 If $r=\infty$, the presentation $B=\QQ[\Scal_1(\bfx),\Scal_2(\bfx),\ldots]$ holds. Similarly define (see \eqref{eq2:Vand})
$$
\sum_{\blamb\in\Pcal_r}\Scal_\blamb(\bfz)s_\blamb(\bfz_r)=\exp(\sum_{i\geq 1}x_ip_i(\bfz_r)),
$$
where $p_i(\bfz_r)=z_1^i+\cdots+z_r^i$ are the Newton power sums.

Since the $\Scal_\blamb(\bfx_r)$ are linearly independent, they form  a $\QQ$-basis of $B_r$:
$$
B_r:=\bigoplus_{\blamb\in\Pcal_r}\QQ\cdot \Scal_\blamb(\bfx_r).
$$

The {\em Hall} inner product on $B$ is defined by  $\langle \Scal_\blamb(\bfx),\Scal_\bmu(\bfx)\rangle=\delta_{\blamb,\bmu}$.
In particular,
\[
\left\langle
\exp\left(\sum_{i\geq 1}x_i p_i(\bfz_r)\right),\,
\exp\left(\sum_{i\geq 1}x_i p_i(\bfw_r)\right)
\right\rangle
=
\sum_{\blamb\in \Pcal_r}s_\blamb(\bfz_r)s_\blamb(\bfw_r),
\]
which is one of the possible spelling  of the classical Cauchy identity as in \cite[Section 4.3, formula (3)]{Fulyoung}. 

\claim{\bf Proposition} \cite[Lemma A.9.4]{Ful}\label{prop:prop24}
{\em Pieri's rule for Schur polynomials holds:
\[
S_i(\bfx)\cdot S_\blamb(\bfx)=\sum_{\bmu\in PF_i(\blamb)} S_\bmu(\bfx).
\]
}
\proof It is enough to put $f(z)=\exp(\sum_{i\geq 1}x_iz^i)=\sum_{i\geq 0}\Scal_i(\bfx)z^i$, and then applying Proposition~\ref{propo:genpier}.\qed

\claim{} If $r=\infty$, then 
$
B:=B_\infty:=\QQ[\bfx]
$ 
is a polynomial ring in the infinitely many indeterminates $\bfx:=(x_1,x_2,\ldots)$ or $E:=(e_1,e_2,\ldots)$, and is often referred to as the {\em bosonic Fock space}.
Let 
$$
\sum_{i\geq 1}\Scal_i(\tilde{\d})t^i:=\exp\left(\sum_{i\geq 1}{t^i\over i}\d_i\right), \qquad \d_i:={\d \over \d x_i}.
$$ 
The equality defines differential operators $\Scal_i(\tilde{\d}):=[t^i]\exp\left(\displaystyle{\sum_{i\geq 1}}\displaystyle{1\over i}\d_it^i\right)$.
If  $\bft_r:=(t_1,\ldots,t_r)$, let $p_i(\bft_r)$ denote, as usual, the corresponding power sum $t_1^i+\cdots+t_r^i$ of degree $i$. More generally, one defines differential operators $\Scal_\blamb(\tilde{\d})$ imposing the equality 
\[
\sum_{\blamb\in \Pcal_r}\Scal_\blamb(\tilde{\d})s_\blamb(\bft_r)=\exp\left(\sum_{i\geq 1}{p_i(\bft_r)\over i}{\d\over \d x_i}\right),
\] 
from which, invoking once more Lemma~\ref{lem:giapie}, $
\Scal_\blamb(\tilde{\d})=\det(\Scal_{\lambda_j-j+i}(\tilde{\d}))$: itis the differential 
operator  obtained by replacing each $x_i$ occurring in $\Scal_\blamb(\bfx)$ by $
\displaystyle{{1\over i}{\d\over \d x_i}}$.

\claim{} To recover Schubert Calculus introduced in the first part of this section, one works on  those subspaces of the ring $B$
\[
\widetilde{B}_{r,n}
:=
\bigoplus_{\blamb\in\Pcal_{r,n}}\QQ\cdot \Scal_\blamb(\bfx),
\]
 which are spanned by Schur basis elements associated to partitions $\blamb\in\Pcal_{r,n}$. This is obviously linearly isomporphic to $B_{r,n}$. The unique linear map
$
\pi_{r,n}:\widetilde{B}_{r,n}\longrightarrow H_*(G(r,n),\QQ)$ sending
$\Scal_\blamb(\bfx)\longmapsto \Omega^\blamb
$
is obviously a vector space isomorphism.

\medskip

\bclm{\bf Theorem.}\label{thm:thm13} {\em The following equality holds in the homology module $H_*(G(r,n),\QQ)$ of the Grassmannian $G(r,n)$:
\begin{equation}
\sum_{i\geq 0 }\sigma_it^i\cap \Omega^\blamb=c_t(\Qcal_r)\cap \Omega^\blamb
=\pi_{r,n}\left(\exp\left(
\sum_{i\geq 1}\frac{t^i}{i}\frac{\partial}{\partial x_i}
\right)\Scal_\blamb(x)\right).
\end{equation}
}
\eclm

\proof We use the fact that, by~\cite[p.~76]{MacDonald},  $\Scal_{\blamb}(\tilde{\d})$ is the  Hall adjoint of $\Scal_\blamb(\bfx)$. Thus:
$$
\delta_{\blamb,\bmu}=\langle \Scal_\blamb(\bfx), \Scal_\bmu(\bfx)\rangle=\langle 1, \Scal_\blamb(\tilde{\d})\Scal_\bmu(\bfx)\rangle=\left. \Scal_\blamb(\tilde{\d})\right|_{x_i=0}\Scal_\bmu(\bfx)
$$
We have to show that
$$
[t^i](c_t(\Qcal_r)\cap \Omega^\blamb)=[t^i]\pi_{r,n}(\exp\left(\sum_{i\geq 1}{t^i\over i}\d_i\right)\Scal_\blamb(\bfx)),
$$
i.e.
that
$$
\sigma_i\cap \Omega^\blamb=\pi_{r,n}\left(\Scal_i(\tilde{\d})S_\blamb(\bfx)\right)
$$
As the left hand side can be computed applying \eqref{eq5:dualp}, the equality follows once one shows that the evaluation $\Scal_i(\tilde{\d})\Scal_\blamb(\bfx)$ obeys to a Pieri's-like formula.
Indeed
\begin{eqnarray*}
\Scal_i({\tilde{\d}})\Scal_\blamb(\bfx)&=&\sum_{\bmu\in\Pcal_r}\big\langle\Scal_i({\tilde{\d}})\Scal_\blamb(\bfx),\Scal_\bmu(\bfx)\big\rangle\Scal_\bmu(\bfx)=\sum_{\bmu\in\Pcal_r}\big\langle \Scal_\blamb(\bfx), \Scal_i(\bfx)\Scal_\bmu(\bfx)\big\rangle\Scal_\bmu(\bfx)\\
&=&\sum_{\bmu\in\Pcal_r}\big\langle \Scal_\blamb(\bfx), \sum_{\bnu\in PF_i(\bmu)}\Scal_\bnu(\bfx)\big\rangle\,\Scal_\bmu(\bfx)\cr\cr
&=&\sum_{\bmu\in\Pcal_r\,|\, \blamb\in PF_i(\bmu)}\Scal_\bmu(\bfx)=\sum_{\bmu\in PF_{-i}(\blamb)}\Scal_\bmu(\bfx)
\end{eqnarray*}
 
having used  Pieri's formulas \ref{prop:prop24} and Remark~\ref{rmk:remark39}.
\qed

\medskip

\claim{\bf Remark.} Notice that the operator $\exp\left(\sum_{i\geq 1}{t^i\over i}\d_i\right)$  is half part, so to speak, of a vertex operator as in \cite[Proposition 5.1]{KaRai}, up to replacing $t^i$ by $t^{-i}$. In addition, Theorem \ref{thm:thm13} shows that Schubert calculus on the Grassmannian is governed by explicit differential operators on the polynomial ring $B$. In particular, the operator $\partial/\partial x_1$ corresponds to the  intersection of the Pl\"ucker image of a Schubert variety with the hyperplane class.

\bclm{\bf Corollary.}\label{cor:cor31} {\em
Let $\blamb\in\Pcal_{r,n}$. Then
\begin{equation}
\sigma_1 \cap \Omega^\blamb
=
\pi_{r,n}\left(\frac{\partial \Scal_\blamb(\bfx)}{\partial x_1}\right).
\end{equation}
In particular, if $d=|\blamb|$, one obtains that the Pl\"ucker degree of the Schubert variety $\Omega^\blamb(F_\bullet)$ is given by
\begin{equation}
\omega_\blamb
=
\frac{\partial^d}{\partial x_1^d}\Scal_\blamb(\bfx).
\end{equation}
which also computes, by \cite[Theorem 2.31]{smirnov}, the number of standard Young tableaux of shape $\blamb$.}
\eclm

\proof
The first identity follows by taking the linear term in $t$ in Theorem~\ref{thm:thm13}. Iterating $d$ times the operator $\partial/\partial x_1$, yields the second identity. Since $\pi_{r,n}$ acts trivially on constants, it can be omitted in the second formula.
\qed

The Pl\"ucker degree of a Schubert variety is a particular case of what one should mean by {\em integral} on a Grassmannian. The cohomology ring is direct sums of homogeneous pieces indexed by (real) codimension:
$$
H^*(G(r,n),\QQ):=\bigoplus_{i=0}^{r(n-r)}H^{2i}(G(r,n),\QQ)
$$
The homogeneous classes in $H^{2r(n-r)}(G(r,n),\QQ)$ will be said top-codimensional classes. The full cohomology $H^*(G(r,n), \QQ)$ is a ring with respect the cup product, which is the linear extension of bilinear maps
$$
H^{2i}(G(r,n),\QQ)\otimes_\QQ H^{2j}(G(r,n),\QQ)\stackrel{\cup}{\lra} H^{2(i+j)}(G(r,n),\QQ)
$$
Each element of the cohomology is a finite sum of homogeneous elements
$$
\alpha_0+\alpha_1+\cdots+\alpha_{r(n-r)}
$$
with $\alpha_i\in H^{2i}(G(r,n), \QQ)$

\bclm{\bf Definition.}
Let \(G(r,n)\) be a Grassmannian and let
\[
\alpha=\alpha_0+\alpha_1+\cdots+\alpha_{r(n-r)}\in H^{*}(G(r,n),\QQ)
\]
Its integral is defined by
\[
\int_{G(r,n)}\alpha\cap \Omega^0=\int_{G(r,n)}\alpha_{r(n-r)}\cap \Omega^0
:=
\deg(\alpha_{r(n-r)}\cap \Omega^0)
\]
which is the degree of $0$-dimensional cycle (roughly speaking formal sum of reduced points with non-negative multiplicities) being the Poincar\'e dual of the the cohomology class $\alpha_{r(n-r)}$ (see e.g. \cite{Ful})
\eclm
Since any cohomology class is a finite linear combination of special Schubert cycles, one can also define integral on a Grasmanniann the evaluation of any  expression
\be
\int_{G(r,n)}\sigma_1^{i_1}\cdots\sigma_k^{i_k}\cap \Omega^{(n-r)^r}
\ee
which is not zero if and only if $i_1+2i_2+\cdots+ki_k=r(n-r)$, which is clearly a generalization of the degree $\int_{G(r,n)}\sigma_1^{r(n-r)}\cap \Omega^0$.
More generally, if $\Omega$ is any homology classe of dimension $d$, one calls integral on the Grassmannian the evaluation
$$
\int_{G(r,n)}\sigma_{1}^{i_1}\cdots\sigma_{k}^{i_k}\cap \Omega
$$
provided $i_1+2i_2+\cdots+ki_k=d$. Therefore the Pl\"ucker degree of a Schubert variety is  a special example of integral on the Grassmanniann $G(r,n)$.
\claim{} \label{sec:sec313} Theorem~\ref{thm:thm13} establishes an important vocabulary, namely that each special Schubert cycle $\sigma_i$ can be understood as a linear differential operators, a polynomial in $\tilde{\d}$ with rational coefficients. This implies that there are more fundamental cohomology classes in $G(r,n)$, not having  the immediate geometrical meaning of the special Schubert cycles. One can define $d_i\in H^*(G(r,n),\QQ)$  via the equality
$$
d_i\cap \Omega^\blamb=\pi_{r,n}\left({\d\over \d x_i}\Scal_\blamb(\bfx)\right)
$$
The coordinate vectorfields ${\d\over \d x_i}$ can be then regarded as cohomology classes, and then the most elementary integral on a Grassmannian are those of the form
$$
\int_{G(r,n)}d_1^{i_1}\cdots d_k^{i_k}\cap \Omega^\blamb=\left\langle 1,\,\,\, {\d^{i_1+\cdots+ki_k}\over \d x_2^{i-2}\cdots\d x_k^{i_k}}\Scal_\blamb(\bfx)\right\rangle
$$
\claim{\bf Remark.} Although the $d_i$s have not any immediate geometric interpretation, we must notice that $d_1$ do in fact coincide with $\sigma_1$. In addition, over a projective space $\PP^n:=G(2, 1+n)$ they represent in fact the dual of the classes of linear subvarieties of codimension $i$. In other words, {\em B\'ezout} theorem for Grassmannian can be expressed through the compact formula
$$
\exp(\sum_{i\geq 1}\d_it^i)\exp(\sum_{i\geq 1}x_ip_i(\bfz_r))=\exp(\sum_{i\geq 1}(t^ip_i(\bfz_r))
$$
meaning that the 
$$
\sigma_i\cap \Omega^\blamb=\pi_{r,n}\left([t^is_\blamb(\bfz_r]\exp(\sum_{i\geq 1}(t^ip_i(\bfz_r))\right).
$$

\claim{} It turns out that the Pl\"ucker degree of a Schubert variety can be 
expressed in determinantal form, leading to a formula that Stanley used in 
\cite[Theorem 7.2.1]{stanley2} to count the number of standard Young 
tableaux of shape $\blamb$. It coincides  with the degree of a Schubert 
variety and also with the dimension of the {\em Specht} module ${\mathbb S}
^\blamb$, i.e. the dimension of the representation of the symmetric group 
$S_d$ associated to the partition $\blamb$ of weight $d$ (See also 
Section~\ref{sec:sec4}). Our key tool is the multiple iteration of a partial 
derivative.

\claim{\bf Proposition.} (Notation as in Lemma \ref{lem:giapie}) Let $d=|\blamb|$. The degree $\omega_\blamb$ of the Schubert variety $\Omega_\blamb(F_\bullet)$ is given by
\be
\omega_\blamb:=d!\cdot \Delta_\blamb(\exp(t))\label{eq:forsoma}
\ee
\proof  Lemma \ref{lem:giapie} and the fact that  
$\exp(\sum_{i\geq 1}x_ip_i(\bfz_r))=\displaystyle{\prod_{i=1}^r}\exp(\sum_{j\geq 1}x_jz_i^j)$  easily implies the Jacobi-Trudi formula
\be
\Scal_\blamb(\bfx_r):=\Delta_\blamb(\Scal_i(\bfx_r))=\det(\Scal_{\blamb_j-j+i}(\bfx_r))_{1\leq i,j\leq r}=\begin{vmatrix}
\Scal_{\lambda_1}(\bfx)&\cdots&\Scal_{\lambda_r+r-1}(\bfx)\cr\vdots&\ddots&\vdots\cr
\Scal_{\lambda_1+r-1}(\bfx)&\cdots&\Scal_{\lambda_r}(\bfx)
\end{vmatrix}.\label{eq2:JTfor}
\ee 
 
Each polynomial $\Scal_i(\bfx)$ occurring in the determinant \eqref{eq2:JTfor} is of the form
\[
\Scal_i(\bfx)=\frac{x_1^i}{i!}+g_i(\bfx),
\]
where $g_i(\bfx)$ is a polynomial in which $x_1$ appears with degree strictly smaller than $i$.
Substituting into (\ref{eq2:JTfor}), the determinant becomes

\begin{equation}
\left|
\begin{matrix}
\displaystyle{x_1^{\lambda_1}\over\lambda_1!}+g_{\lambda_1} &
\displaystyle{x_1^{\lambda_2-1}\over (\lambda_2-1)!}+g_{\lambda_2-1} &
\cdots &
\displaystyle{x_1^{\lambda_r+r-1}\over (\lambda_r+r-1)!}+g_{\lambda_r+r-1}\\[6pt]
\displaystyle{x_1^{\lambda_1+1}\over (\lambda_1+1)!}+g_{\lambda_1+1} &
\displaystyle{x_1^{\lambda_2}\over \lambda_2!}+g_{\lambda_2} &
\cdots &
\displaystyle{x_1^{\lambda_r+r-2}\over (\lambda_r+r-2)!}+g_{\lambda_r+r-2}\\
\vdots & \vdots & \ddots & \vdots\\
\displaystyle{x_1^{\lambda_1+r-1}\over (\lambda_1+r-1)!}+g_{\lambda_1+r-1} &
\displaystyle{x_1^{\lambda_2+r-2}\over (\lambda_2+r-2)!}+g_{\lambda_2+r-2} &
\cdots &
\displaystyle{x_1^{\lambda_r}\over \lambda_r!}+g_{\lambda_r}
\end{matrix}
\right|.
\label{eq:hugedet}
\end{equation}
Let $d=|\blamb|$. A few  manipulations show that (\ref{eq:hugedet}) can be written in the form
\[
x_1^d\Delta_\blamb(\exp(t)) + F(\bfx),
\]
where $F(\bfx)$ is a polynomial in which $x_1$ occurs with degree strictly smaller than $d$. By Corollary \ref{cor:cor31}
\[
\omega_\blamb
=
\frac{\partial^d}{\partial x_1^d}\left(x_1^d\Delta_\blamb(\exp(t)) + F(\bfx)\right)
=
d!\cdot \Delta_\blamb(\exp(t)).
\]
as claimed.\qed
%
%
%
%
%

\medskip
\bclm{\bf Example.}
The degree of the Schubert variety indexed by $(3,2,1)$ is
\[
\omega_{(3,2,1)}=16,
\]
which coincides with the number of standard Young tableaux of shape\,\, $\tiny{\bmtx\ydiagram{3,2,1}\cr\emtx}$.
\eclm

\claim{} Degrees of Schubert varieties obeys two kind of recursion. The former is 
\be
\omega_\blamb={\d^d\over \d x_1^d} \Scal_\blamb(\bfx) ={\d^{d-1}\over \d x_1^{d-1}}\sum_{i=1}^(X^r(\blamb-\ep_j)=\sum_{j=1}^r \omega_{\blamb-\ep_j}\label{eq:FulFor}
\ee
where by $\ep_j=(0,\ldots,1,\ldots,0)$ we denote the row with $r$-entries all $0$ but $1$ in the $j$th spot, $\blamb-\ep_j$ is the partition with all parts equal to those of $\blamb$ but the $j$-th one which is $\lambda_j-\ep_j$,  setting $\omega_{\blamb-\ep_j}=0$ if $\blamb-\ep_j$ is not a partition. Formula~\eqref{eq:FulFor} is nothing but \cite[Example 4.7.11]{Ful}. For instance in case of partitions of length $2$ we have the following {\em Pascal triangle} for the degree of the corresponding Schubert varieties
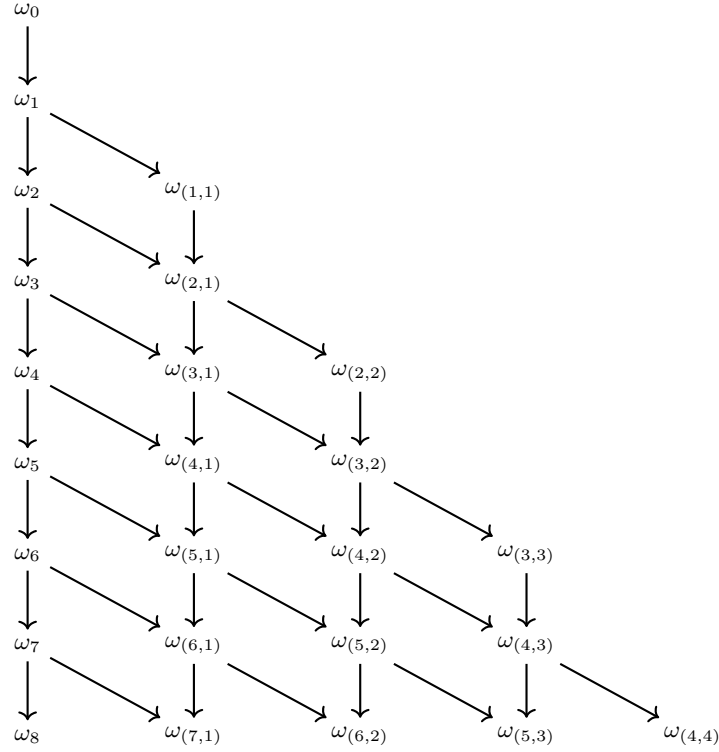
\begin{center}
\begin{tikzpicture}[
  every node/.style={font=\small},
  arr/.style={->, thick}
]

\def\N{8}        
\def\Jmax{8}     
\def\dx{2.2}     
\def\dy{1.2}     

\foreach \i in {0,...,\N} {
  \foreach \j in {0,...,\Jmax} {
    \ifnum\j>\i\relax
    \else
      \pgfmathtruncatemacro{\s}{\i+\j}
      \ifnum\s>8\relax
      \else
        \node (w-\i-\j) at ({\dx*\j},{-\dy*\s})
        {%
          \ifnum\j=0
            $\omega_{\i}$
          \else
            $\omega_{(\i,\j)}$
          \fi
        };
      \fi
    \fi
  }
}

\foreach \i in {0,...,\N} {
  \pgfmathtruncatemacro{\ip}{\i+1}
  \ifnum\ip>\N\relax\else
    \foreach \j in {0,...,\Jmax} {
      \pgfmathtruncatemacro{\s}{\i+\j}
      \pgfmathtruncatemacro{\sp}{\ip+\j}
      \ifnum\j>\i\relax\else
        \ifnum\s>8\relax\else
          \ifnum\sp>8\relax\else
            \draw[arr] (w-\i-\j) -- (w-\ip-\j);
          \fi
        \fi
      \fi
    }
  \fi
}

\foreach \i in {0,...,\N} {
  \foreach \j in {0,...,\Jmax} {
    \pgfmathtruncatemacro{\jp}{\j+1}
    \pgfmathtruncatemacro{\s}{\i+\j}
    \pgfmathtruncatemacro{\sp}{\i+\jp}
    \ifnum\jp>\Jmax\relax\else
      \ifnum\i<\jp\relax\else
        \ifnum\s>8\relax\else
          \ifnum\sp>8\relax\else
            \draw[arr] (w-\i-\j) -- (w-\i-\jp);
          \fi
        \fi
      \fi
    \fi
  }
}

\end{tikzpicture}

-- \small{The ``Pascal triangle'' for the degrees of Schubert varieties associated to partitions of length $2$} --
\end{center}
Recall that $\omega_{n,n}$ is the $n$-th Catalan number $C_n$, Therefore  the above diagram is related with the papers \cite{CataNiede, CataTaise}.

\noindent
The second recursion formula we have alluded to is:
\bclm{\bf Corollary.} Let $\blamb\in\Pcal_r$ such that $|\blamb|=d$. The degree $\omega_\blamb$ of a Schubert variety satifies the following recursive formulas
$$
\omega_\blamb=\sum_{i=1}^r(-1)^{i-1}{d\choose \lambda_{i}+i-1}\omega_{\blamb^{(i)}}
$$
\eclm
\proof
Given the relation
$$
\omega_\blamb=d!\Delta_\blamb(\exp(t))
$$
it is enough to use the expansion of the determinant:
\begin{eqnarray*}
d!\Delta_\blamb(\exp(t))&=&d!\left(\sum_{i=1}^r(-1)^{i-1}{1\over (\lambda_i-i+1)!}\Delta_{\blamb^{(i)}}(\exp(t) \right)\cr\cr
&=&\sum_{i=1}^r(-1)^{i-1}{d \over (\lambda_i-i+1)!|\blamb^{(i)}|!}|\blamb^{(i)}|!\Delta_{\blamb^{(i)}}(\exp(t))=\sum_{i=1}^r{d\choose \lambda_i-i+1}\omega_{\blamb^{(i)}}.\qed
\end{eqnarray*}

\claim{\bf Remark.} In the case $\blamb=(\lambda_1\geq \lambda_2)$
$$
\omega_{(\lambda_1,\lambda_2)}={|\lambda|\choose \lambda_1}-{|\lambda|\choose \lambda_2-1}
$$
In particular, for the partition $(n,n)$:
$$
\omega_{(n,n)}={2n\choose n}-{2n\choose n-1}.
$$
which is the $n$-th Catalan number.

\bclm{\bf Corollaries} (see  \cite[Section 4.3 (4)]{Fulyoung}, \cite[p.~13]{Manivelsym})
{\em

i) For each $r\geq 1$, one has:\qquad
$
\displaystyle{\sum_{d\geq 0}\frac{t^d}{d!}\sum_{\substack{\blamb\in\Pcal_r\\|\blamb|=d}}\omega_\blamb s_\blamb(\bfz_r)}
=
\exp\bigl(t(z_1+\cdots+z_r)\bigr).
$

\medskip
In particular,
$
(z_1+\cdots+z_r)^d
=
\sum_{\blamb\vdash d}\omega_\blamb s_\blamb(\bfz_r).
$
In addition,  ii), 
 $$
\sum_{|\blamb|=d}\omega_\blamb^2=d!
$$

}
\eclm

\proof
We have that
\begin{eqnarray*}
\sum_{|\blamb|=d}\langle x_1^d, \Scal_\blamb(\bfx)\rangle s_\blamb(\bfz_r)&=&\langle x_1^d, \sum_\blamb \Scal_\blamb(\bfx)s_\blamb(\bfz_r)\rangle=\cr\cr&=&{\d^d\over \d x_1^d}_{\bfx=0}\exp(\sum x_ip_i(\bfz_r))=p_1(\bfz_r)^d\cr\cr &=&(z_1+\cdots+z_r)^d
=\sum_{|\blamb|=d}\omega_\blamb s_\blamb(\bfz_r).
\end{eqnarray*}

Therefore\quad
$
\displaystyle{\sum_{d\geq 0}{t^d\over d!}\sum_{|\blamb|=d}\omega_\blamb s_\blamb(\bfx)=\sum_{d\geq 0}{t^d\over d!}(z_1+\cdots+z_r)^d=\exp(t(z_1,\ldots+z_r))}
$.

To see ii), one first expands $x_1^d$ as a linear combination of Schur functions
\be
x_1^d=\sum_{|\blamb|=d}\langle x_1^d, S_\blamb(\bfx)\rangle S_\blamb(\bfx).\label{eq3:x_1d}
\ee
Taking the $d$-th derivative with respect to $x_1$ on either side of the equality \eqref{eq3:x_1d}:

$$
d!={\d^d\over \d x_1^d}x_1^d=\sum_{|\blamb|=d}\big\langle 1, {\d^d\over \d x_1^d}S_\blamb(\bfx)\big\rangle{\d^d\over \d x_1^d}\Scal_\blamb(\bfx)=\sum_{|\blamb|=d}\omega_\blamb^2.
$$

\qed
%

\claim{} More generally, one has
\[
\sigma_i \cap \Omega^\blamb
=
\pi_{r,n}\left(\Scal_i(\widetilde{\d})\Scal_\blamb(\bfx)\right).
\]

For example,
\[
\sigma_3 \cap \Omega^\blamb
=
\pi_{r,n}\left[
\left(
\frac{1}{6}\frac{\partial^3}{\partial x_1^3}
+
\frac{1}{2}\frac{\partial^2}{\partial x_1 \partial x_2}
+
\frac{1}{3}\frac{\partial}{\partial x_3}
\right)
\Scal_\blamb(\bfx)
\right].
\]


It is worth of remarking that the Schur polynomials $\Scal_i(\bfx)$ satisfy the following relation
\be
{\d \Scal_i(\bfx)\over \d x_j}={\d^j \Scal_i(\bfx)\over \d x_1^j}=\Scal_{i-j}(\bfx)\label{eq2:si-j}
\ee
as a straightforward output of the  identity
\begin{eqnarray*}
{\d \Scal_i(\bfx)\over \d x_j}&=&[z^i]{\d \over \d x_j}\left(\sum_{k\geq 0}\Scal_k(\bfx)z^k\right)
=[z^i]{\d\over \d x_j}\exp(\sum x_iz^i)\cr\cr &=&[z^i]\left(z^j\exp(\sum_{i\geq 1}x_kz^k)\right)=[z^i]\left(\sum_{k\geq 0}\Scal_{k+j}(\bfx)z^k\right)=\Scal_{i-j}(\bfx)
\end{eqnarray*}
where in the last equality we have invoked Remark \ref{rmk:remark39}.

\bclm{\bf Example.} Let $G:=G(2,n+2)$ be the Grassmannian of the $2$-planes  in $\CC^{n+2}$.  Its dimension is $2n$ and its fundamental class is $\Omega_{n,n}$ which can be modeled by the Schur polynomial $\Scal_{n,n}(\bfx)$. Therefore, if $a+2b=2n$, the action of $\sigma_1^a\sigma_2^{b}$ against $\Omega^{n,n}$ is a top codimension cycle,  an {\em integral} on the Grassmannian, also written as
$$
\int_{G(2,n+2)}\sigma_1^a\sigma_2^b\cap[G]
$$
Its computation is straightforward using the differential calculus formalism. The translation of $\sigma_1$ and $\sigma_2$ into differential operators is
$$
\sigma_1={\d\over \d x_1}\qquad \mathrm{and }\qquad \sigma_2={1\over 2}\left({\d^2\over \d x_1^2}+{\d\over \d x_2}\right)
$$
The sought for number is then given by
$$
\int_{G}\sigma_1^{2n-2b}\sigma_2^b\cap [G]={1\over 2^b}{\d_1^{2n-2b}}\left({\d_1^2}+{\d_2}\right)^b\Scal_{n,n}={1\over 2^b}\sum_{j=0}^b{b\choose j}\d_1^{2b-2j}\d_2^j\Scal_{n,n}(\bfx)
$$
Writing 
$$
\Scal_{n,n}=\Scal_n^2-\Scal_{n-1}\Scal_{n+1}
$$
$$
\d_2^j(\Scal_n^2-\Scal_{n-1}\Scal_{n+1})=\sum{j\choose k}(\Scal_{n-2k}(\bfx)\Scal_{n-2j+2k}(\bfx)- \Scal_{n+1-2k}(\bfx)\Scal_{n-1+2j+2k}
$$
Now let
$$
\d_1^{2b-2j}(\Scal_{n-2k}(\bfx)\Scal_{n-2j+2k}(\bfx)- \Scal_{n+1-2k}(\bfx)\Scal_{n-1+2j+2k})=\eta_{n-2k, n-2j+2k}
$$
One then obtains
$$
\int_{G}\sigma_1^{2n-2b}\sigma_2^b\cap [G]={1\over 2^b}\sum_{j=0}^b\sum_{k=0}^j{b\choose j}{j\choose k}\eta_{n-2k, n-2j+2k}
$$
where $\eta_{n-2k, n-2j+2k}=\omega_{n-2k, n-2j+2k}$ if $n-2k\geq n-2j+2k$, is $0$ if $n-2k=m-2j+2k-1$ and is $-\omega_{n-2j+2k+1,n-2k+1}$ if $n-2k<m-2j+2k-1$

\eclm

%
%
%
%

\medskip

\medskip


\section{Review on Schubert Derivations}\label{sec:sec4}

\claim{} Doing Schubert Calculus by means of partial derivatives acting on 
$B$ amounts to identify the latter with the universal container of all the 
homologies of the Grassmannians at once.  Then the partial derivatives 
play the role of cohomology operators. Indeed they generate the cohomology 
as a $\QQ$-algebra, because all the classical Schubert (co)cycles, once 
they are identified with the differential operators $\Scal_\blamb(\tilde{\d})$, are 
explicit polynomial expressions in the partial derivatives. The other 
point of view, that we propose here, by further developing ideas borrowed from \cite{SCHSD,ESC,SDIWP}, is letting all the $B_r$ playing the role of the 
cohomologies of $G(r,\infty)$ and studying their action on another 
universal containers of the homologies, which is the exterior algebra of a 
countably finite vector space, once  identified with that of 
polynomials with rational coefficients. To this goal, it is particular relevant  to  work not just with $B$ (although this would be clearly possible) but with 
all the $B_r$ together, for finite $r$. On such rings 
it would not make sense to consider partial derivatives ${\d\over \d x_i}$ 
for $i>r$ (because $B_r$ is generated only by $x_1,\ldots, x_r$) but there exists  
suitable replacement coming from natural endomorphisms of the exterior 
algebra, as will be shown next. These operators enable the construction of an 
algebra that mimics, and is related to, the usual Weyl algebra of polynomial rings.

\claim{\bf Exterior Algebra.} \label{sec:sec34} In the sequel we consider 
$V:=\QQ[X]$ as being the vector space of polynomials with rational 
coefficients and, occasionaly, $V_n:=\displaystyle{\QQ[X]\over (X^n)}$, so 
that $V_\infty:=V$. For $n\in\NN\cup\{\infty\}$, the natural basis of 
$V_n$ is $\bfX:=(X^i)_{0\leq i<n}$. Let ${\bm\d}:=(\d^i)_{i\geq 0}$, where
$$
\d^j:={1\over j!}\left.{d^j\over d X^j}\right|_{X=0}: V\sra \QQ,
$$
which is clearly the unique linear form on $V$ such that $\d^j(X^i)=\delta^j_i$. The vectorspace $V_n^*=\bigoplus_{i\geq 0}
\QQ\d^i$ is the {\em restricted dual} of $V_n$ (the ordinary dual space 
if $n<\infty$).  One has a canonical projection 
$$
\rho_n:V\sra V_n
$$ which in the sequel will also understood as the endomorphism of $V$ obtained by its composition with the canononical injection $V_n\hookrightarrow 
V$. L
The generating function of the basis $\bfX$ and ${\bm\d}$ are given by
$$
\bfX(z)=\sum_{i\geq 0}X^iz^i=\exp\left(\sum_{i\geq 1}X^i{z^i\over i}\right)\qquad \mathrm{and}\qquad \bmd(t):=\sum_{j\geq 0}\d^jt^j=\left.\exp\left(t{d\over d X}\right)\right|_{X=0},
$$ 
The Taylor formula gives $\d(w)\bfX(z)=\displaystyle{\sum_{i\geq 1}}t^iz^i$, as expected.

\claim{} To the vector space $V$ one associates a sequence of {\em exterior powers} 
$(\bw^iV)_{i\geq 0}$, where one sets $\bw^0V=\QQ$, $\bw^1V=V$ and, for $r\geq 2$, $\bw^rV$ is understood  as  the vector space generated by all the expressions 
$X^{i_1}\w\cdots\w X^{i_r}$, modulo the relations putting to zero any 
such element, if there are $j\neq k$ such that $X^{i_j}=X^{i_k}$. 
One will write 
$
\bw^rV=\bigoplus_{\blamb\in\Pcal_r}\QQ\cdot \bfX^r(\blamb)
$,
where
$$
\bfX^r(\blamb):=X^{r-1+\lambda_1}\w X^{r-2+\lambda_2}\w\cdots\w 
X^{\lambda_r}
$$
is the natural basis element of $\bw^rV$ uniquely attached to the partition $
\blamb\in\Pcal_r$.
The {\em exterior algebra} of $V$ is
$$
\bw V=\left(\bigoplus_{r\geq 0}\bw^rV,\w\right),
$$
where $\w$ is the juxtaposition product. The {\em Hall scalar product} on $\bw  V$ is defined by $\langle \bfX^r(\blamb), \bfX^s(\bmu)\rangle=\delta_{\blamb,\bmu}$.

\claim{} \label{secs:sec44} Exterior algebras are strictly 
related to  symmetric polynomials. The wedge product\linebreak
$
\bfX(z_1)\w\cdots\w \bfX(z_r)
$
is evidently a linear combination of the vectors $\bfX^r(\blamb)$. Its obvious 
divisibility by the Vandermonde determinant
$
\Delta_0(\bfz_r):=\prod_{i<j}(z_i-z_j)
$ defines the Schur symmetric polynomials $s_\blamb(\bfz_r)$
through the equality
\be
\bfX(z_1)\w\cdots\w \bfX(z_r)
=
\Delta_0(\bfz_r)\sum_{\blamb\in \Pcal_r}\bfX^r(\blamb)s_\blamb(\bfz_r).\label{eq:for_pre_schur}
\ee

\bclm{\bf Definition.} {\em A {\em Hasse-Schmidt derivation} on $\bw V$ is a $\QQ$-linear map $\Dcal(z)=\sum_{i\geq 0}D_iz^i:\bw V\sra \bw V\llb z\rrb$ such that $D_0=\id_{\bw V}$ and
\be
\Dcal(z)(u\w v)=\Dcal(z)u\w \Dcal(z)v.
\ee
}
\eclm
By equating the formal expansions of $\Dcal(z)(u\w v)=	\sum_{i\geq 0}D_i(u\w v)z^i$ with that of $\Dcal(z)u\w \Dcal(z)v=\sum_{i\geq 0}D_iu\cdot z^i\w \sum_{j\geq 0}D_jv\cdot z^j$, one easily sees that $\Dcal(z)=\sum_{i\geq 0}D_iz^i\in\End_{\QQ}(\bw V)$ is a $HS$-derivation on $\bw V$ if and only if, for all $i\geq 0$,
$
D_i(u\w v)=\sum_{j=0}^iD_ju\w D_{i-j}v.
$
The condition $D_0=\id_{\bw V_n}=\id_{\bw V_n}$ implies that   $\Dcal(z)$ is invertible as $\End_\QQ(\bw V)$-valued formal power series. We denote by $\ovDcal(z)$  the inverse of $\Dcal(z)$ as $\End_\QQ(\bw V)$-valued formal power series,  and write it as
\be 
\ovDcal(z)=\sum_{i\geq 0} (-1)^i\ovD_iz^i
\ee
 Motivated by \cite{GSCH}, we call the {\em Cayley-Hamilton} relations
the obvious equalities
\begin{eqnarray}
\ovDcal(z)(\Dcal(z)u\w v)&=&u\w \ovDcal(z)v\label{eq:CH1}\\
\Dcal(z)(u\w \ovDcal(z))&=&\Dcal(z)u\w v.\label{eq:CH2}
\end{eqnarray}

\claim{} According to the notation of Lemma~\ref{lem:giapie}, we denote by $\Dcal(\bfz_r)$ the product $\Dcal(z_1)\cdots\Dcal(z_r)$. It is a $\w$-algebra homomorphism as well, in the sense that $\Dcal(\bfz_r)(u\w v)=\Dcal(\bfz_r)u\w \Dcal(\bfz_r)v$.
Clearly $\Dcal(\bfz_r)$ is an $\End_\QQ(\bw V)$-valued formal power series in symmetric polynomials which admits the same kind of expansion as in Lemma~\ref{lem:giapie}:
$
\Dcal(\bfz_r)=\sum_{\blamb\in\Pcal_r}\Delta_\blamb(\Dcal)s_\blamb(\bfz_r)
$,
which intrinsecally define endomorphisms $\Delta_\blamb(\Dcal)$ of $\bw V$. By Lemma \ref{lem:giapie},  $\Delta_\blamb(\Dcal)$ is an explicit determinant in the coefficients of $\Dcal(z)$:
\be 
\Delta_\blamb(\Dcal)= \det(\Dcal_{\lambda_j-j+i})
\ee
and the purely combinatorial Pieri's formula holds  by invoking once more Proposition~\ref{propo:genpier}
\be
D_i\Delta_\blamb(\Dcal)=\sum_{\bmu\in PF_i}\Delta_\bmu(\Dcal).
\ee
 It is not surprising that $HS$-derivations are related to derivations of the exterior algebra. We stick  ourselves to the easiest case of {\em even} derivations,  shortly said {\em derivations}, accoding to the following:
\bclm{\bf Definition.} \label{def47:derivation} {{\em  A derivation} on the exterior algebra $\bw V$ is a $\QQ$-linear map
$d:\bw V\sra \bw V$ such that $d(u\w v)=du\w v+u\w dv$.}
\eclm 
Let $\bfd(z):=\sum_{i\geq 1}d_iz^i:\bw V\sra \bw V\llb z\rrb$. Then
$$
\bfd(z)(u\w v)=\bfd(z)u\w v+u\w \bfd(z)v\quad \iff \quad d_i(u\w v)=d_iu\w v+u\w d_iv,\qquad \forall i>0 \quad \forall u,v\in\bw V
$$
in which case, by abuse of terminology, we say  that $\bfd(z)$ is a derivation as well. 

If the coefficients $d_i$ of $\bfd(z)$ are pairwise commuting, i.e. $[d_i,d_j]=0$ for all $i,j\geq 1$, then one may speak of its exponential
\be 
\exp(\bfd(z))=1 + \sum_{j\geq 1}{\bfd(z)^j\over j!}:\bwV\sra \bwV\llb z\rrb
\ee
\bclm{\bf Proposition.}\label{prop6: expd}
{\em The map $\Dcal(z)=\exp(\bfd(z))\in \End_\QQ(\bw V)$ is a $HS$-derivation on $\bw V$ if and only if 
$\bfd(z)$ is a derivation.}
\eclm
\proof Suppose that $\bfd(z)$ is a derivation. Then $\Dcal(z):=\exp(\bfd(z))=1+\bfd(z)+\displaystyle{\bfd(z)^2\over 2!}+\cdots$ and, in general, $\bfd(z)^k(u\w v)=\sum_{j=0}^k{k\choose j}\bfd(z)^ju\w \bfd(z)^{k-j}v$. Therefore
\begin{eqnarray*}
\Dcal(z)(u\w v)&=&\sum_{n\geq 0}{1\over n!}\sum_{k=0}^n \bfd(z)^k(u\w v)=\sum_{n\geq 0}{1\over n!}\sum_{k=0}^n{n\choose k}\bfd(z)^ku\w \bfd(z)^{n-k}v\cr\cr
&=&\sum_{n\geq 0}\sum_{h+k=n}{1\over  k!}\bfd(z)^ku\w {1\over h!}\bfd(z)^{h}v=\exp(\bfd(z)u\w\exp(\bfd(z))v.
\end{eqnarray*}
Conversely, suppose $\Dcal(z)$ is a $HS$- derivation. We contend that its logarithm is a derivation.
We have
\[
{d\over dz}\left[\Dcal(z)(u\w v)\right]={d\over dz}\Big(\Dcal(z)u\w \Dcal(z)v\Big).
\]
Therefore
\[
\Dcal(z)'(u\w v)= \bfd'(z)\Dcal(z)(u\w v)=\bfd'(z)\Dcal(z)u\w \Dcal(z)v+\Dcal(z)u\w \bfd(z)'\Dcal(z)v
\]
It is now sufficient to multiply either side of the above equality by $\ovDcal(z)$ to obtain:
$$
\bfd(z)'(u\w v)=\bfd(z)'u\w v+u\w \bfd(z)'v,
$$
which proves that $\bfd'(z)=\sum_{i\geq 1}{id_i}z^i$ is a derivation, i.e. each $d_i$ is a derivation and  $\bfd(z)$ is a derivation.\qed.
\claim{} \label{sec:sec49} Let $gl(V)$ be the Lie sub-algebra of $\End_Q(V)$ of all the endomorphisms of $V$ such that $\d^j(\phi(X^i))=0$ if $|i-j|\gg0$. Identifying $gl(V)$ with an algebra of matrices, this is the same as the Lie algebras of matrices with finitely many non zero diagonals \cite[p.~34]{KaRai}.  Clearly $gl(V_n)=\End_\QQ(V_n)$.

\bclm{\bf Definition.} {\em 
If $\phi\in gl(V)$, we denote by $\updelta(\phi)$ the unique derivation of the exterior algebra such that  $\delta(\phi)u=\phi(u)$ for all $u\in V$. }
\eclm

\claim{\bf Remark.} An easy check shows that the map $\updelta:gl(V_n)\sra \End(\bw V_n)$ is a Lie algebra homomorphism, in the sense that $\updelta([\phi, \psi])=[\updelta(\phi),\updelta(\psi)]$.

\claim{}\label{sec:dadj} Because trivially
$
[\updelta(\phi^i),\updelta(\phi^j)]=0
$
in $\End(V)$, by Proposition~\ref{prop6: expd}, 
$$
\exp(\Dcal^\phi(z))=\exp\left(\sum_{j\geq 1}{\delta(\phi^j)\over j}z^j\right):\bw V\sra \bw V\llb z\rrb
$$
is the unique  HS-derivation  such that $\Dcal^\phi(z)_{|V}=\sum_{n\geq 0}\phi^n z^n$.
We denote by  $\phi^*$ its $\langle\, , \rangle$--adjoint with respect to the Hall scalar product, i.e. 
$
\langle \phi^*(u), v\rangle=\langle u, \phi(v)\rangle
$
A standard exercises, based on induction on the degree of the exterior power, shows that  adjoint of a derivation of $\bw V$ is a derivation and $\delta(\phi^*)$ is a derivation as well. As $\delta(\phi^*)$ and $\delta(\phi)^*$ agree on $V=\bw^1V\subseteq \bw V$, they do coincide.

\claim{\bf Example.} As in Section~\ref{sec:sec34}, let us think of  $\rho_n:V\sra V$ thought of as the composition of the projection $V\sra V_n$ and the canonical injection $V_n\subseteq V$, which identifies $V_n$ with the vector subspaces generated by $(X^0, X^1,\ldots, X^{n-1})$. Clearly $\ker (\rho_n)=\Span_\QQ(X^i)_{i\geq n}$. The derivation $\updelta(\rho_n)$ is defined as follows
$$
\updelta(\rho_n)(u\w v)=\updelta(\rho_n)u\w v+u\w \updelta(\rho_n)v
$$
and $\delta(\rho_n)\bfX^r(\blamb)=0$ if and only if $X^{r-i+\lambda_i}\notin V_n$.

For example
$
\updelta(\rho_2)(X^5\w X^2)=X^5\w X^2
$
whereas $\updelta(\rho_2)(X^5\w X^3)=\updelta(\rho_2)X^5\w X^3+X^5\w \delta(\rho_2)X^2=0$. Clearly, if $n=\infty$, $\updelta(\rho_\infty)$ is the identity. 
\medskip
 
Let us now turn to a special distinguished case. Let $X:V\sra V$ defined by $u\mapsto Xu$.
\bclm{\bf Definition.}
{\em The Schubert derivations $\sigma_+(z),\ovsig_+(z):\bw V\sra \bw V\llb z\rrb $ are defined by
\begin{eqnarray*}
\sigma_+(z)&=&\Dcal^{X}(z)=\sum_{i\geq 0}\sigma_iz^i=\exp\left(\sum_{i\geq 1}{\delta(X^i)\over i}z^i\right)\cr\cr
\ovsig_+(z)&=&\ovDcal^{X}(z)=\sum_{i\geq 0}\sigma_iz^i=\exp\left(-\sum_{i\geq 1}{\delta(X^i)\over i}z^i\right)
\end{eqnarray*}
}
\eclm
The derivations $\sigma_+(z)$ and $\ovsig_+(z)$ are clearly mutually inverses. The map $\sigma_+(z)$ is the unique HS-derivation such that 
$$
\sigma_+(z)u=\sum_{i\geq 0}X^iu\cdot z^i
$$
while $\ovsig_+(z)$ is the unique such that
$$
\ovsig_+(z)u=u-Xu\cdot z.
$$
\bclm{\bf Proposition.} {\em The following equality holds in $\bw^rV_n\llb \bfz_r\rrb$
\be 
\sigma_+(\bfz_r)\bfX^r(0)=\sum_{\blamb\in\Pcal_r}
\bfX^r(\blamb)s_\blamb(\bfz_r)=\sum_{\blamb\in\Pcal_{r,n}}\bfX^r(\blamb)s_\blamb(\bfz_r)\label{eq3:fakVa}
\ee
}
\eclm
\proof Notation as in formula~\eqref{eq:for_pre_schur}. For $r=1$, we set $\Delta_0(z_1)=1$, so that formula \eqref{eq3:fakVa} holds in this case. For $r=2$ we have:
\begin{eqnarray*}
\sigma_+(z_1)X^0\w \sigma_+(z_2)X^0&=&\sigma_+(z_1,z_2)\big(\ovsig_+(z_2)X^0\w\ovsig_+(z_1)X^0\big)\cr\cr
&=&\sigma_+(z_1,z_2)\big[(X^0-z_1X^1)\w (X^0-z_2X^1)\big]\cr\cr
&=&(z_1-z_2)\sigma_+(z_1,z_2)X^1\w X^0
\end{eqnarray*}
and thus the property holds for $r=2$. Then we argue by induction. Let us suppose that 
$$
\sigma_+(z_1)X^0\w\cdots\w\sigma_+(z_{r-1})X^0=\Delta_0(\bfz_{r-1})\sigma_+(\bfz_r)\bfX^{r-1}(0).
$$
Then, on one hand $\sigma_+(z_1)X^0\w\cdots\w\sigma_+(z_{r-1})X^0\w\sigma_+(z_r)X^0$ is clearly a multiple of $\prod_{i<j}(z_i-z_j)$. On the other hand, using induction:
\begin{eqnarray*}
\sigma_+(z_1)X^0\w\cdots\w\sigma_+(z_{r-1})X^0\w\sigma_+(z_r)X^0=\Delta_0(\bfz_{r-1})\sigma_+(\bfz_{r-1})\bfX^{r-1}(0)\w \sigma_+(z_r)X^0
\end{eqnarray*}
Invoking the Cayley-Hamilton relations \eqref{eq:CH1} we get, from the last side:
\begin{eqnarray}
&&\Delta_0(\bfz_{r-1})\sigma_+(\bfz_{r-1})\bfX^{r-1}(0)\w \sigma_+(z_r)X^0\cr\cr
&=&\Delta_0(\bfz_{r-1})\sigma_+(\bfz_r)\ovsig_+(z_r)\bfX^{r-1}(0)\w \ovsig_+(\bfz_r\hskip-2pt\setminus z_r)X^0\cr\cr
&=&\sum_{i=0}^{r-1}(-1)^it^i\bfX^{r-1}(1^i)\w \big(X^0-e_1(\bfz_{r-1}X^1+\cdots+(-1)^{r-1}e_{r-1}(\bfz_{r-1})X^{r-1}\big).\cr&&\label{eq3:redvndm}
\end{eqnarray}
Expression~\eqref{eq3:redvndm} is clearly a multiple of $\bfX^r(0)$, which is divisible by $\prod_{i=1}^r (z_i-z_r)$. To determine its coefficient one looks at that of $z_1\cdots z_{r-1}=e_{r-1}(\bfz_{r-1})$ in \eqref{eq3:redvndm}. In all cases one gets a summand of the form
\be
(-1)^{r-1}z_1\cdots z_{r-1} \bfX^{r-1}(0)\w X^{r-1}.\label{eq3:exch}
\ee
Now:  
$$
\bfX^{r-1}(0)\w X^{r-1}=(-1)^{r-1} X^{r-1}\w\bfX^{r-1}(0)=(-1)^{r-1}\bfX^r(0).
$$
 Therefore
$$
(-1)^{r-1}e_{r-1}(\bfz_{r-1})=(-1)^{r-1}z_1\cdots z_{r-1} \bfX^{r-1}(0)\w X^{r-1}=z_1\cdots z_{r-1}\bfX^r(0).
$$

proving that \eqref{eq3:redvndm} is equal to $\Delta_0(\bfz_r)\bfX^{r}(0)$, as claimed.\qed

\claim{} Recall notation  from  Section \ref{sec:sec41}.  One defines 
$$
\exp(\sum_{i\geq 1}x_iz^i)u=\sigma_+(z)u,\qquad \forall u\in \bw V.
$$
Therefore each $x_i$ can be regarded as a derivation of the exterior algebra
$$
x_i(u\w v)=x_iu\w v+u\w x_i v
$$
for all $r\geq 1$.

\bclm{\bf Proposition.}
{\em The map $B_r\sra \bw^rV$ mapping $x_i\mapsto \displaystyle{\updelta(X^i)\over i}$ makes $\bw^rV$ into a free $B_r$-module of rank $1$, generated by $\bfX^r(0)$, such that $\Scal_\blamb(\bfx)\bfX^r(0)=\bfX^r(\blamb)$.}
\eclm
\proof
One proves it at the level of generating functions.
\begin{eqnarray*}
\sum_{\blamb\in \Pcal_r}\Scal_\blamb(\bfx_r)\bfX^r(0)s_\blamb(\bfz_r)&=&\exp\left(\sum_{i\geq 1}x_iz^i\right)\bfX^r(0)=\prod_{j=1}^r\exp(\sum_{i\geq 1}x_1z_j^i)\bfX^r(0)\cr
&=&\sigma_+(\bfz_r)\bfX^r(0)=\sum_{\blamb\in\Pcal_r}\Delta_\blamb(\sigma)\bfX^r(0)=\sum_{\blamb\in\Pcal_r}\bfX^r(\blamb)
\end{eqnarray*}

\claim{\bf Remark.} It is obvious that there must exist one, and hence 
infinitely many, linear isomorphisms $B_r\sra \bw^rV$. The reason why we 
preferred to phrase the distinguished one as above, is to model Poincar\'e 
duality for Grassmannian. The map $B_r\sra \bw^rV$ can be understood as a 
finite type version of the Boson-Fermion correspondence. Classically this is 
an isomorphism between the ring $B$ and a fixed charge vector subspace of an 
infinite wedge power. Morally, one can think of it as an isomorphism between 
$B$ and the vectorspace generated by the symbols $\bfX^\infty(\blamb)$, one 
for each partition, to be understood as the projective limit in the category 
of graded vector spaces of the exterior power. In other words, let $
(\bw^rV)_w=\bigoplus_{\blamb\vdash w}\QQ \bfX^r(\blamb)$. If $r\gg w$, then 
$(\bw^{r}V)_w$ and $(\bw^{r+s}V)_w$ have the same dimension and the map
$(\bw^{r+s}V)_w\sra (\bw^{r}V)_w$ given by $\bfX^{r+s}(\blamb)\mapsto 
\bfX^r(\blamb)$ is the required isomorphism. One then considers $(\bw^\infty 
V)_w=\lim_{\sra}(\bw^rV)_w$, with basis $(\bfX^\infty(\blamb))_{\blamb\vdash 
w}$. One sets $\bw^\infty V=\bigoplus_{w\geq 0}(\bw^\infty V)_w$. The vacuum 
vector is $\bfX^\infty(0)$ which is approximated by the finite 
approximations $\bfX^r(0)$. The Hall inner product is $\langle 
\bfX^\infty(\blamb),\bfX^\infty(\bmu)\rangle=\delta_{\blamb,\bmu}$, and if 
$u\in \bw^\infty V$, then
$$
\langle u \rangle:=\int_{\bw^rV}u:=\langle \bfX^\infty(0), u\rangle
$$
is said to be the {\em vacuum expectation value} of the vector $u$ and, as a matter of fact, can be thought as an integral on the infinite Grassmannian.  By abuse of terminology, we can say the coefficient of $\bfX^r(0)$ in an element $u\in \bw^rV$ is its vaccum expectation value.
One can also set the following
\bclm{\bf Definition.}
{\em The equality
\be
\int_{\bw^\infty}\sigma_1^{i_1}\cdots\sigma_k^{i_k}=\left\langle \sigma_{-1}^{i_1+2i_2+\cdots+ki_k}\bfX^\infty(0)\right\rangle \label{eq:fermint}
\ee
will be said integral on the infinite Grassmannian.
}
\eclm
By the boson-fermion correspondence, evaluating \eqref{eq:fermint}  is the same as computing
$$
{\d^d\over \d x_1^d}\Scal_1(\bfx)^{i_1}\cdots\Scal_{i_k}(\bfx)^{i_k}
$$

Going back to geometry, for finite $r$,
the {\em vacuum vector} $\bfX^r(0)$ plays the role of the fundamental 
class of the Grassmannian $G(r,\infty)$, and $\Scal_\blamb(\bfx_r)\cong 
\Delta_\blamb(\sigma_+)$ play the role of cohomologies classes: the action 
of $B_r$ on $\bw^rV$ plays the role of the {\em cap} product $\cap: 
H^*(G(r,n),\QQ)\otimes_\QQ H_*(G(r,n),\QQ)\sra H_*(G(r,n),\QQ)$.

Again, as a consequence of Lemma~\ref{lem:giapie}, one obtains one more Pieri's-like formula, i.e.
\be 
 \sigma_i\bfX^r(\blamb)=\sum_{\bmu\in PF_i}\bfX^r(\bmu)
 \ee
Indeed $\sigma_i\bfX^r(\blamb)=\sigma_i\Delta_\blamb(\sigma)\bfX^r(0)$ and then one observes that $\Delta_\blamb(\sigma)$ is the coefficient of the Schur polynomial $s_\blamb(\bfz_r)$ in the expansion of $\sigma_+(\bfz_r)$.

\claim{} Let $X^{-i}\in \End(V)$ be the adjoint of $X^i$ with respect to the Hall scalar product, i.e.
$$ 
\langle X^{-i}u,v\rangle= \langle u,X^iv\rangle \qquad \forall u,v\in V.
$$
with respect the Hall inner product on $\bw V$.
It is a locally nilpotent endomorphism whihc cna be explicitly described by noticing that $X^{-i}X^k=X^{k-i}$ if $k\geq i$ and $0$ otherwise.
Then, by~\ref{sec:dadj}, $\delta(X)^*=\delta(X^{-i})$. Therefore $\delta(X^{-i})$ is the unique derivation of the exterior algebra such that $\delta(X^{-i})u=X^{-i}u$. One can consider the associated HS-derivations, i.e
$$
\sigma_-(t)=\sum_{j\geq 0}\sigma_jt^j:=\exp\left(\sum_{i\geq 1}{\delta(X^{-i})\over i}t^i\right)
$$
as well as  its inverse
$$
\ovsig_-(t)=\sum_{j\geq 0}(-1)^j\ovsig_jt^j=\exp\left(-\sum_{i\geq 1}{\delta(X^{-i})\over i}t^i\right).
$$
They are the unique HS-derivations of $\bw V$ such that $\sigma_-(t)X^j=\sum_{i=1}^j{X^{j-i}t^i}$ and $\ovsig_-(t)X^j=X^j-X^{j-1}t$.
Recall that $B_r$ is generated as a $\QQ$-algebra by $\Scal_i(\bfx_r)$. Let us define $\delta(X^{-i})$ on $B_r$ as follows
$$
(\delta(X^{-i})p)\bfX^r(0)=\delta(X^{-i})(p\bfX^r(0))
$$
So, for instance $(\updelta(X^{-j})\Scal_i(\bfx_r))\bfX^r(0)=\updelta(X^{-j})\bfX^r((i))=\delta(X^{-j})X^{r-1+i}\w \bfX^{r-1}(0)$.

\bclm{\bf Lemma.} {\em For all $i,j\geq 0$ and all $r\geq 1$, one has
$$
\delta(X^{-i})\Scal_j(\bfx_r)=\Scal_{j-i}(\bfx_r).
$$
}
\eclm
\proof  First of all we check that for all $r\geq 1$, $\delta(X^{-i})\bfX^r(0)=0$ if $i\geq 1$. Indeed 
$$
\delta(X^{-i})\bfX^r(0)=\sum_{\blamb\in\Pcal_r}\Big\langle \delta(X^{-i})\bfX^r(0), \bfX^r(\blamb)\Big\rangle\bfX^r(\blamb)=\Big\langle \bfX^r(0), \delta(X^i)\bfX^r(\blamb)\Big\rangle\bfX^r(\blamb)=0
$$
the last vanishing due to the fact that for all $\blamb\in\Pcal_r$, $\delta(X^i)\bfX^r(\blamb)$ is a linear combination of $\bfX^r(\bmu)$ with $|\bmu|>0$. As a consequence
\begin{eqnarray*}
(\delta(X^{-i})\Scal_j(\bfx_r))\bfX^r(0)&=&\delta(X^{-i})X^{r-1+j}\w \bfX^{r-2}(0)\cr\cr
&=&(\delta(X^{-i})X^{r-1+j})\w \bfX^{r-2}(0)+ X^{r-1+j}\w \delta(X^{-i}\bfX^{r-2}(0))\cr\cr
&=&\bfX^{r-1+j-i}\w \bfX^{r-2}(0)=\Scal_{j-i}(\bfx_r)\bfX^r(0)
\end{eqnarray*}

which proves the claim.\qed

 Let us now denote by $J:=(j_1,\ldots, j_r)$ a multi-index of non-negative integers of length $r$ and by $\Scal_J(\bfx_r)=\prod_{k=1}^r\Scal_{j_k}(\bfx_r)$.
Let $\ep_k$ denotes the coordinate row $(\delta_{ik})_{1\leq i\leq r}$, all entries $0$ but $1$ in position $j$.

\bclm{\bf Proposition.} {\em The derivation  $\delta(X^{-i})$ of the exterior algebra $\bw V$ enjoys the Leibniz rule when evaluated against products of the form $\Scal_{J}(\bfx_r)\in B_r$, i.e.
$$
\delta(X^{-i})\big(\Scal_{j_1}(\bfx_r)\cdots\Scal_{j_r}(\bfx_r)\big)=\sum_{k=1}^r\Scal_{J-i\ep_k}(\bfx_r),
$$
}

\eclm 
\proof
First of all we notice that if $\blamb\in\Pcal_r$, one has
$$
\delta(X^{-i})\Scal_\blamb(\bfx_r)=\sum_{k=1}^r\Scal_{\blamb-i\ep_k}(\bfx_r)
$$
where $\Scal_{\blamb-i\bfe_k}(\bfx_r)=0$ if $\blamb-i\ep_k$ is not a partition. Let $|J|=j_1+\cdots+j_r$. An easy check based on applying the Leibniz rule enjoyed by the derivative of a determinant with respect to its columns,  shows that
$$
\Scal_J(\bfx_r)=\sum_{|\blamb|=|J|}\Scal_\blamb(\bfx_r),
$$
whence
\begin{eqnarray*}
\delta(X^{-i})\Scal_J(\bfx_r)&=&\sum_{|\blamb|= |J|}\delta(X^{-i})\Scal_\blamb(\bfx_r)
=\sum_{|\blamb|= |J|}\sum_{k=1}^r\Scal_{\blamb-i\ep_k}(\bfx_r)\cr\cr
&=&\sum_{k=1}^r \sum_{|\blamb|= |J|}\Scal_{\blamb-i\ep_k}(\bfx_r)=\sum_{k=1}^r\Scal_{J-i\ep_k},
\end{eqnarray*}
which proves the claim.\qed
\bclm{\bf Remark.} It must be emphasized that $\delta(X^{-i})$ satisfies Leibiniz's rule only if applied to the product of no more than $r$ non unital algebra generators of $B_r$. For instance, in $B_2$, 
$$
\delta(X^{-1})\Scal_1^3(\bfx_2)=\delta(X^{-1})(\Scal_2(\bfx_2)\Scal_{1}(\bfx_2))=\Scal_1(\bfx_2)^2+\Scal_2(\bfx_2)\neq 3\Scal_1(\bfx_2)^2.
$$ 

\eclm

\bclm{\bf Proposition.} {\em If $r=\infty$, then $\delta(X^{-i}):B\sra B$ is a $\QQ$-derivation, coinciding with the partial derivative $\displaystyle{\d\over \d x_i}$.}

\eclm
\proof Indeed, $B$ is generated as a $\QQ$-algebra by $\Scal_j(\bfx)$, and $\delta(X^{-i})$ satisfies Leibniz rules against the product of arbitrary numbers of generatorrs. Therefore $\delta(X^{-i})$ is a derivation. Since
$$
{\d\over \d x_i}\exp(\sum_{j\geq 1}x_jz^j)=z^i\exp(\sum_{j\geq 1}x_jz^j)
$$
it follows that 
$$
{\d\Scal_i(\bfx)\over \d x_j}=\Scal_{j-i}(\bfx)=\delta(X^{-i})\Scal_j(\bfx)
$$
therefore
$$
\delta(X^{-i})={\d \over \d x_i}
$$
as the two operators coincide when evaluated against the generators $\Scal_j(\bfx)$ of the algebra $B$. \qed

We notice that while the partial derivatives   $\d_i={\d\over \d x_i}$ are 
defined in $B$ for all $i\geq 1$, the operator $\d_j$ makes no sense in In 
$B_r$ if $j>r$. Therefore $\delta(X^{-i})$ are the correct replacement of 
the partial derivatives in $B_r$. It follows that $B_r$ is acted on by the 
creations operators ``multiplication'' by $x_i$ as well as by the 
annihilation operators $\delta(X^{-i})$. It is natural to wonder if such 
operators define a Weyl agebra as they do in $B$.
We see that the operators $\delta(X^{-i})$ and $x_j$ acts on $\bw^rV_n$.
Since the trace homomorphism $\updelta$ is a Lie-algebra homomorphism, one trivially has the commutation relations:
\be
[\updelta(X^i),\updelta(X^j)]=\delta([X^i,X^j]).
\ee
 Notice that $\delta(X^0)$ acts as the identity in $\bw V$ and  is 
therefore a central element in the sense that $\updelta(X^0)
\updelta(X^j)=\updelta(X^j)\updelta(X^0)$. Moreover $\updelta(X^i)$ acts on 
$\bw V$ as annihilation operators if $i<0$ and creation operators if $i>0$.
We have that $\updelta(X^i)$ is a derivation on $\bw V$ for all $i\in \ZZ$. A few easy experiments show that the operators
$
\delta(X^i)
$, $i\in\ZZ$,
do not satisfy the usual Heisenberg relations on the whole exterior algebra. The failure is measured by explicit kind of {\em boundary operators}.

\bclm{\bf Proposition.} \label{prop:prop35} {\em Let $V=\QQ[X]$. If  $0\leq i\leq j$, the sub-algebra $W^\wedge$ of $
\End_\QQ(\bw^rV)$ generated by $\delta(X^i)$ for all $i\in\ZZ$, satisies the 
following commutation relations:
$$
[\updelta(X^i),\updelta(X^j)]=\updelta(X^{j-i}\rho_{i-1}).
$$
}
\eclm
\proof 
The proof amounts to a simple case by case analysis. We will prove that for all $i,j\geq 0$ and  all $k\geq 0$ the equality
$
[X^{-i}, X^j]X^k=X^{j-i}\rho_{i-1}(X^k)
$
holds.
One has
$$
[X^{-i}, X^j]X^k=X^{-i}X^{j+k}-X^jX^{-i}X^k
$$
\begin{enumerate}
\item[{\bf case 1.}] If $k\geq i$, then $\rho_{i-1}(X^k)=0$. In this case one has
$$
X^{-i}X^{j+k}-X^jX^{-i}X^k=X^{j+k-i}-X^{j+k-i}=0=X^{j-i}\updelta(\rho_{i-1}(X^k)
$$
\item[{\bf case 2.}] Let $k\leq i-1$. Then$$
X^{-i}X^{j+k}-X^jX^{-i}X^k=X^{j+k-i}=X^{j-i}X^k=X^{j-i}\rho_{i-1}(X^k)
$$
where $X^{j-i}X^k$ is zero if $j-i+k<0$.
\end{enumerate}
Therefore we have proven, for all $k\geq 0$ and all $i,j>0$, the equality
$
[X^{-i},X^j]X^k=X^{j-i}\rho_{i-1}(X^k)
$
which in turns implies the equality
$
[X^{-i},X^j]=X^{j-i}\rho_{i-1}
$
holding in $gl(V_n)$. Thus
$$
\updelta([X^{-i},X^j])=[\updelta(X^{-i}, X^j]=\updelta(X^{j-i}\rho_{i-1})
$$
and the Proposition is proven. \qed

\bclm{\bf Corollary.} {\em Let $m,n$ be positive. If $r\geq max(m,n)$, then
$$
[\updelta(X^{-m}),\updelta(X^n)]\bfX^r(0)=m\delta_{m,n}\bfX^r(0).
$$
}
\eclm
\proof
Under the hypotheses of the statement, by Proposition~\ref{prop:prop35}
\begin{eqnarray}
[\updelta(X^{-m}),\updelta(X^n)]\bfX^r(0)&=&\updelta(X^{n-m}\rho_{-m-1})\bfX^r(0)\cr\cr
&=&X^{r-1}\w\cdots\w X^{m}\w \updelta(X^{n-m}\rho_{-m-1})\bfX^{m-1}(0)\cr \nonumber\\
&=&X^{r-1}\w\cdots\w X^{m}\w \sum_{j=0}^{m-1}X^{m-1}\w\cdots\w X^{j+n-m}\w\cdots X^0.
\end{eqnarray}
 If $n>m$, since $r\geq n$ one has
$
j+n-m\leq (m-1)+n-m=n-1\leq r-1
$
which means that $X^{j+n-m}$ is an exterior factor of $\bfX^r(0)$ for all $0\leq j\leq m-1$. Therefore the expression is zero.
If $n<m$ instead we have two cases. Either $j+n-m<0$, and then the corresponding summand is zero ($X^{n-m}X^j=0$ if $j<n-m$). Or $0\leq j+n-m\leq m-1$, hence $X^{j+n-m}$ is an exterior factor of $\bfX^m(0)$ and then the corresponding summand vanishes.
We are finally left with the case $m=n$. The commutation relations read in this case as
$$
[\updelta(X^{-m}), \delta(X^m)]\bfX^r(0)=\delta(\rho_{m-1})\bfX^r(0).
$$
Using Leibniz rule, $\delta(\rho_{m-1})$ acts as the identity to each factor in between $0$ and $m-1$. So, it counts how many $X^j$ with $0\leq j\leq m-1$ do occur in $\bfX^r(0)$. These  are exactly $m$. This prove the Corollary.\qed

\bclm{\bf Corollary.} Let $\blamb\in \Pcal_r$, $m,n>0$ and  $s\gg (r, |m|, |n|)$. Then
$$
[\updelta(X^{m}), \updelta(X^{n})]\bfX^s(\blamb)=-m\delta_{m,-n}\bfX^r(\blamb)
$$
\eclm
\proof The proposition is obvious if $m,n$ are both positive or negative, as in that case the operators $\delta(X^m)$ and $\delta(X^n)$ commute. If $m<0$ and $n>0$ let us write $m=-p$, with $p$ positive. Then
\begin{eqnarray*}
[\updelta(X^{-p}), \updelta(X^{n})]\bfX^s(\blamb)&=&\updelta(X^{n-p}\rho_{p-1})\bfX^s(\blamb)\cr\cr
&=&X^{s-1+\lambda_1}\w\cdots\w X^{s-r+\lambda^r}\w\updelta(X^{n-p}\rho_{p-1})\bfX^{s-r-1}(0)\cr\cr
&=&X^{s-1+\lambda_1}\w\cdots\w X^{s-r+\lambda_r}\w p\bfX^{s-r-1}(0)\cr\cr
&=&p\bfX^s(\blamb)
\end{eqnarray*}

\bclm{\bf Remark.} {If $r=\infty$ one obtains
$$
[\updelta(X^m),\updelta(X^n)]\bfX^\infty(\blamb)=-m\delta_{m,-n}\bfX^\infty(\blamb)
$$
Using the boson-fermion correspondence, the fact that $\delta(X^{-i})={\d\over \d x_i}$ on $B$, that $\delta(X^i)$ correspond to the multiplication by $ix_i$, one recover the classial Heisenberg algebra relations as in \cite[p.~12]{KaRai} 
}
\eclm

The previous result shows that the classical Heisenberg algebra appears as the stable limit of the finite-rank commutation relations. The operators
\[
\delta(X^{-m})
\qquad\text{and}\qquad
\delta(X^m)
\]
should therefore be regarded as annihilation and creation operators respectively, while the projections $\rho_{m-1}$ encode the boundary corrections responsible for the failure of the Heisenberg relations at finite rank.
\claim{\bf Connection to the representation theory of the symmetric group}. In this final part we point out a connection between derivations on exterior algebras, in the sense of Definition~\ref{def47:derivation}, 
and the representation theory of the symmetric group.

Let \(S_d\) be the symmetric group on \(d\) letters.As well known, to every partition
$\blamb\vdash d$  one attach the {\em Specht module} \(\mathbb S^\blamb\),  which turns out to be an 
irreducible representation of $S_d$. The partition $\blamb$ prescribes the 
cycle type of a partition. The Specth module $\bbS^\blamb$ admits a basis of 
so-called standard polytabloids, and there is one polytabloid for each 
standard Young tabloid. A standard Young tabloid of shape $\blamb\vdash d$ is a filling of the 
Young diagram of $\blamb$ with the integers $1,2,\ldots,n$ each appearing exactly once, such 
that entries increase from left to right along rows and from top to bottom along columns, and 
is known that 
$
\dim \mathbb \bbS^\blamb=\omega_\blamb,
$
where \(\omega_\blamb\), the degree of the Schubert variety $\Omega^\blamb(F_\bullet)$ as in Section~\ref{sec3:diff}, also coincides with the number of standard Young tableaux of shape \(\blamb\), which can be computed by the hook-length formula.

By definition of representation, any \(\tau\in S_d\) acts linearly on \(\mathbb S^\lambda\). It will then be viewed as an endomorphism of $\bbS^\blamb$, whose  trace  is, by definition, the value of the irreducible character associated with $\blamb$:
\[
\chi^\blamb(\tau)
=
\operatorname{Tr}\left(\tau_{|\bbS^\lambda}\right).
\]
By \cite[Chapter 7]{Ful} that characters are constant on conjugacy classes, this value depends only on the
cycle type of \(\tau\) which is prescribed by the chosen partition. If
$
\bmu=(1^{m_1}2^{m_2}\cdots)\vdash d
$
is the corresponding cycle type, it is customary to denote such a  value by \(\chi^\blamb_\bmu\). 
\bclm{\bf Theorem.} \label{thm:thm431}{\em Let
$
\bmu=(1^{m_1}2^{m_2}\cdots)\vdash d
$
and let \(\blamb\vdash d\) be a partition of length at most \(r\). Then
\[
\delta(X^{-1})^{m_1}\delta(X^{-2})^{m_2}\cdots
\bfX^r(\blamb)
=
\chi^\blamb_\bmu\,\bfX^r(0),
\]
where \(\chi^\blamb_\bmu\) denotes the value of the irreducible character of \(S_d\)
indexed by \(\blamb\), evaluated on the conjugacy class of cycle type \(\bmu\).
}
\eclm

\proof
Recall, from Section \ref{secs:sec44} that
\[
\sigma_+(\bfz_r)\bfX^r(0)=
\sum_{\lambda\in\Pcal_r}
\bfX^r(\blamb)s_\blamb(\bfz_r)
=
\frac{\bfX(z_1)\wedge\cdots\wedge \bfX(z_r)}
{\Delta_0(\bfz_r)}.
\]

If $i\geq 1$,
$
\delta(X^{-i})\bfX(z)
=
\sum_{m\geq i}X^{m-i}z^m
=
z^i\bfX(z),
$
from which, invoking the Leibniz rule enjoyed by each $\delta(X^{-i})$,
\[
\delta(X^{-i})\sigma_+(\bfz_r)\bfX^r(0)
=
p_i(\bfz_r)\sigma_+(\bfz_r)\bfX^r(0),
\]

We now consider an exterior algebra  analogous of the radial vector field:
\[
\rho(\bfx_r)
=
\exp\Bigl(\sum_{i\geq 1}x_i\delta(X^{-i})\Bigr).
\]

One clearly has
$
\rho(\bfx_r)\sigma_+(\bfz_r)\bfX^r(0)
=
\exp\Bigl(\sum_{i\geq 1}x_ip_i(\bfz_r)\Bigr)
\bfX^r(0).
$
Expanding both sides,
\[
\sum_{\blamb\in\Pcal_r}
\rho(\bfx)\bfX^r(\blamb)\,
s_\blamb(\bfz_r)
=
\exp\Bigl(\sum_{i\geq 1}x_ip_i(\bfz_r)\Bigr)
\sum_{\blamb\in\Pcal_r}
\bfX^r(\blamb)s_\blamb(\bfz_r).
\]
One then takes the {\em vacuum expectation value}, i.e. the coefficient of the vacuum vector \(\bfX^r(0)\), obtaining
\be
\langle \bfX^r(0), \sum_{\blamb\in\Pcal_r}
\,\rho(\bfx_r)\bfX^r(\blamb)\,
s_\blamb(\bfz_r)\rangle
=
\exp\Bigl(\sum_{i\geq 1}x_ip_i(\bfz_r)\Bigr)=\sum_{\blamb\in\Pcal_r}\Scal_\blamb(\bfx_r)s_\blamb(\bfz_r)\label{eq4:32}
\ee
Now we give a look to the full expansion of the left hand side.
Since the derivations \(\delta(X^{-i})\) commute with each other, we have
\[
\rho(\bfx)
=
\exp\left(\sum_{i\geq 1}x_i\delta(X^{-i})\right)
=
\prod_{i\geq 1}\exp\left(x_i\delta(X^{-i})\right),
\]
whence
\[
\rho(\bfx)
=
\prod_{i\geq 1}
\left(
\sum_{m_i\geq 0}
\frac{x_i^{m_i}}{m_i!}\delta(X^{-i})^{m_i}
\right)=\sum_{\mu=(1^{m_1}2^{m_2}\cdots)}
\left(
\prod_{i\geq 1}\frac{x_i^{m_i}}{m_i!}
\right)
\delta(X^{-1})^{m_1}
\delta(X^{-2})^{m_2}\cdots .
\]

As for the right hand side, we invoke the Frobenius' character formula as displayed e.g. in  \cite[Section 6, formula (12)]{Fulyoung} recalling that the $p_i$ there occurring are related to our $x_i$ by the equality $x_i=p_i/i$, obtaining:
\be
\Scal_\blamb(\bfx)=\sum_{\mu=(1^{m_1}2^{m_2}\cdots)}
\chi^\blamb_\bmu
\prod_{i\geq 1}
\frac{x_i^{m_i}}{m_i!}\label{eq4:eq33}
\ee
By replacing \eqref{eq4:eq33} the expression in formula \eqref{eq4:32}:
\[
\exp\Bigl(\sum_{i\geq 1}x_ip_i(\bfz_r)\Bigr)
=\sum_{\blamb}\Scal_\blamb(\bfx)s_\blamb(\bfz_r)=
\sum_{\blamb}
\left(
\sum_{\mu=(1^{m_1}2^{m_2}\cdots)}
\chi^\blamb_\bmu
\prod_{i\geq 1}
\frac{x_i^{m_i}}{m_i!}
\right)
s_\blamb(\bfz_r).
\]

whence
\be
\sum_{\mu=(1^{m_1}2^{m_2}\cdots)}
\left(
\prod_{i\geq 1}\frac{x_i^{m_i}}{m_i!}
\right)
\delta(X^{-1})^{m_1}
\delta(X^{-2})^{m_2}\cdots=\sum_{\blamb}
\left(
\sum_{\mu=(1^{m_1}2^{m_2}\cdots)}
\chi^\blamb_\bmu
\prod_{i\geq 1}
\frac{x_i^{m_i}}{m_i!}
\right)
s_\blamb(\bfz_r)\label{eq4:finalcom}
\ee
It just remains to compare the coefficients of 
$
\prod_{i\geq 1}\frac{x_i^{m_i}}{m_i!}\,
s_\blamb(\bfz_r),
$
on both sides of~\eqref{eq4:finalcom}, to get:
\[
\big\langle
\delta(X^{-1})^{m_1}
\delta(X^{-2})^{m_2}
\cdots
\bfX^r(\lambda), \bfX(0)\rangle
=
\chi^\blamb_\bmu.
\]

Since our hypthesis was \(|\bmu|=|\blamb|\), the operator
\[
\delta(X^{-1})^{m_1}
\delta(X^{-2})^{m_2}
\cdots
\]
lowers the weight of $\bfX^r(\lambda)$ by exactly $|\blamb|$. Hence the resulting vector
has weight $0$, and therefore is a scalar multiple of \(\bfX^r(0)\). The previous equality
shows that the scalar is precisely \(\chi^\blamb_\bmu\). Thus
\[
\delta(X^{-1})^{m_1}
\delta(X^{-2})^{m_2}
\cdots
\bfX^r(\lambda)
=
\chi^\blamb_\bmu\,\bfX^r(0),
\]
as claimed.\qed

\bclm{\bf Remark.}
Theorem \ref{thm:thm431} shows that the vacuum expectation
\[
\langle
\delta(X^{-1})^{m_1}
\delta(X^{-2})^{m_2}\cdots
\bfX^r(\lambda),
\bfX^r(0)
\rangle
\]
plays the role of a fermionic integral. Thus, irreducible character values of symmetric
groups arise as fermionic analogues of the ordinary intersection integrals on
Grassmannians.
\eclm

\bclm{\bf Example.} Let us compute the value of the character $\chi^{(3,2)}_{(4,1)}$. By the theorem this is just
\[
\left\langle \bfX^2(0),\,\, \delta(X^{-4})\delta(X^{-1})\bfX(3,2)\right\rangle
\]
Now
$$
\delta(X^{-4})\delta(X^{-1})\bfX(3,2)=\delta(X^{-4})\delta(X^{-1})X^4\w X^2=\delta(X^{-1})X^0\w X^2=X^0\w X^1=-X^1\w X^0
$$
from which $\chi^{(3,2)}_{(4,1)}=-1$. 
\eclm

This method illustrated above is rather fast compared to the definition which should construct first the Specht module $\bbS^{(3,2)}$ and then a  $5\times 5$ matrix to determine the trace. However one can also argue that one may evaluate
$$
{\d^2\over \d x_1 \d x_4}\Scal_{(3,2)}(\bfx)
$$ 
where one should use the Jacobi-Trudi formula to express $\Scal_{(3,2)}(\bfx)$ as $\Scal_3(\bfx)\Scal_2(\bfx)-\Scal_1(\bfx)\Scal_4(\bfx)$ and then using Leibniz rules and equalities.

The  method we have presented represents an elementary algebraic alternative to the purely combinatorial celebrated Murnagham-Nakayama formula, as e.g. in \cite[p.~117 and p.~134, Notes and References]{MacDonald}.

\bibliographystyle{amsplain}

\newpage
\bigskip
\noindent
{\rm André~Contiero}\\
{\tt \href{mailto:contiero@mat.ufmg.br}{contiero@mat.ufmg.br}}\\
{\it Universidade Federal de Minas Gerais}\\
{\it Belo Horizonte, MG, Brazil}\

\bigskip
\noindent
{\rm Letterio~Gatto}\\
{\tt \href{mailto:letterio.gatto@polito.it}{letterio.gatto@polito.it}}\\
{\it Dipartimento~di~Scienze~Matematiche}\\
{\it Politecnico di Torino}\

\bigskip
\noindent
{\rm Parham~Salehyan}\\
{\tt \href{mailto:psalehyan@ibilce.unesp.br}{psalehyan@ibilce.unesp.br}}\\
{\it Departamento de Matemática -- IBILCE--UNESP}\\
{\it Campus de S.~J.~ Rio Preto, São Paulo, Brazil}\\
{\it Politecnico di Torino}\
\end{document}